\documentclass[hidelinks,onefignum,onetabnum]{siamart220329}
\usepackage[T1]{fontenc}
\usepackage[utf8]{inputenc}
\usepackage[english]{babel}
\usepackage{mathtools}
\usepackage{algorithm}
\usepackage[noEnd=true]{algpseudocodex}
\usepackage{accents}
\usepackage{multicol}
\usepackage{enumitem}
\usepackage{url}
\usepackage{cleveref}
\usepackage{siunitx}
\usepackage{caption}
\usepackage{color}
\usepackage{setspace}

\newcommand{\runinhead}[1]{\vspace{1mm}\noindent\textbf{#1}\hspace{0.5ex}}
\newcommand{\hS}[1]{\hspace{#1pt}}

\newcommand{\R}{\mathbb{R}}
\newcommand{\N}{\mathbb{N}}

\newcommand{\dx}{\,\mathrm{d}x}
\newcommand{\ds}{\,\mathrm{d}s}
\newcommand{\dt}{\,\mathrm{d}t}
\newcommand{\dmu}{\,\mathrm{d}\mu}
\newcommand{\mc}{\text{mc}}
\newcommand{\card}[1]{\# #1}
\newcommand{\abs}[1]{|#1|}
\newcommand{\Expec}[1]{\operatorname{E}[#1]}
\newcommand{\Var}[1]{\operatorname{Var}[#1]}

\newcommand{\Params}{\mathcal{P}}
\newcommand{\Dt}{{\Delta t}}
\newcommand{\ROM}{\ensuremath{\textsc{rom}}}
\newcommand{\RB}{\ensuremath{\textsc{rb}}}
\newcommand{\ML}{\ensuremath{\textsc{ml}}}

\newcommand{\overcirc}[1]{\accentset{\circ}{#1}}

\newcommand{\code}[1]{\textnormal{\texttt{\detokenize{#1}}}}

\newcommand{\grayout}[1]{{\color{gray}#1}}
\usepackage{dcolumn}
\newcolumntype{s}{D{,}{\cdot}{-1}}
\newcommand{\sci}[2]{#1,10^{#2}}
\newcommand{\scigray}[2]{\grayout{#1},\grayout{10^{#2}}}
\newcommand{\mcol}[2]{\multicolumn{1}{#1}{#2}}

\newsiamthm{example}{Example}
\newsiamthm{assumption}{Assumption}

\usepackage{lipsum}
\usepackage{amsfonts}
\usepackage{graphicx}
\usepackage{epstopdf}

\ifpdf
	\DeclareGraphicsExtensions{.eps,.pdf,.png,.jpg}
\else
	\DeclareGraphicsExtensions{.eps}
\fi

\crefname{section}{Section}{Sections}
\crefname{subsection}{Subsection}{Subsections}

\setitemize{leftmargin=15pt}

\headers{A certified adaptive RB-ML-ROM surrogate for PPDEs}{B.~Haasdonk, H.~Kleikamp, M.~Ohlberger, F.~Schindler, T.~Wenzel}

\title{A new certified hierarchical and adaptive RB-ML-ROM surrogate model for parametrized PDEs%
\thanks{{Funded by BMBF under contracts 05M20PMA and 05M20VSA.
Funded by the Deutsche Forschungsgemeinschaft (DFG, German Research Foundation) under Germany’s Excellence Strategy EXC 2044 - 390685587, Mathematics M\"unster: Dynamics – Geometry – Structure, and EXC 2075 - 390740016.
We acknowledge the support by the Stuttgart Center for Simulation Science (SimTech).
We acknowledge initial experiments on certified VKOGA models by Nicole Neis.}}}

\makeatletter
\newcommand{\specificthanks}[1]{\@fnsymbol{#1}}
\makeatother

\author{Bernard~Haasdonk\thanks{Institute for Applied Analysis and Numerical Simulation, University of Stuttgart, Germany.}
	\and Hendrik~Kleikamp\thanks{Mathematics Münster, University of Münster, Germany.}
	\and Mario~Ohlberger\footnotemark[3]
	\and Felix~Schindler\textsuperscript{\specificthanks{3},}\thanks{Corresponding author: \email{felix.schindler@wwu.de}}
	\and Tizian~Wenzel\footnotemark[2]}

\ifpdf
\hypersetup{
	pdftitle={A new certified hierarchical and adaptive RB-ML-ROM surrogate model for parametrized PDEs},
	pdfauthor={B. Haasdonk, H. Kleikamp, M. Ohlberger, F. Schindler, T. Wenzel}
}
\fi

\begin{document}

\maketitle

\begin{abstract}
  We present a new surrogate modeling technique for efficient approximation of input-output maps governed by param\-etrized PDEs.
  The model is hierarchical as it is built on a full order model (FOM), reduced order model (ROM) and machine-learning (ML) model chain.
  The model is adaptive in the sense that the ROM and ML model are adapted on-the-fly during a sequence of parametric requests to the model.
  To allow for a certification of the model hierarchy, as well as to control the adaptation process, we employ rigorous a posteriori error estimates for the ROM and ML models.
  In particular, we provide an example of an ML-based model that allows for rigorous analytical quality statements.
  We demonstrate the efficiency of the modeling chain on a Monte Carlo and a parameter-optimization example.
  Here, the ROM is instantiated by Reduced Basis Methods and the ML model is given by a neural network or a VKOGA kernel model.
\end{abstract}



\section{Introduction}
\label{sec:intro}

We are concerned with a new efficient and certified ``on demand'' learning approach to approximate parametrized
input-output maps, where the intermediate state is governed by time-dependent parametrized partial differential equations (PDEs).
Approximation of such input-output maps can be addressed by projection-based model reduction
in combination with rigorous a posteriori error control.
Projection-based model reduction for parametrized systems is a very active research area that has seen enormous development in the
past two decades. For an overview we refer to the collections and handbooks
\cite{MR3701994,MR3672144,zbMATH07204874,zbMATH07204875,zbMATH07204877}.
Traditionally, global model reduction approaches aim at approximating the solution over the whole input space in a suitable
offline-stage where the model is generated. In the subsequent online-stage the reduced model can then be queried for any
parameter in an online-efficient way.
However, such global models are typically expensive to construct as they require many snapshots, i.e.
solution instances at different parameter choices.
Also, such global models may be oversized for the problem at hand. For example, in the case of parameter
optimization, a certain trajectory of parameters will evolve, which will not cover the whole parameter domain.

Therefore, adaptive enrichment in model order reduction (MOR) has been suggested to incrementally refine initially coarse ROMs.
Adaptive enrichment has first been proposed in the context of localized MOR \cite{MR3431132}, where successively also
concepts to combine offline training and online enrichment were discussed \cite{MR4014784,MR3702345,zbMATH07358914,Zahr2015,La2014}.
In the context of
particular multi-query applications, such as PDE constrained parameter optimization, adaptive enrichment has recently seen great success
also for global model reduction approaches, e.g.~in combination with trust-region optimization \cite{MR3716566,MR4269464,BKMOSV22}.

In parallel to such efficient model-based reduction methods, purely data-based approaches for machine learning of surrogate models have
increasingly been developed and mathematically investigated in recent years, see e.g.~
\cite{MR4268857,Fresca2022,UP2021,WKH2014,KSHH2018,SH16b,wenzel2021analysis,KKLOO22}.
Despite data driven (supervised) approaches, there is also an increasing interest in unsupervised learning approaches for the solution of PDEs.
In particular, the physics-informed neural networks (PINNs) \cite{MR3881695,MR4229292} received a lot of attention with applications to a large variety of
PDEs. Very recently, an error analysis for PINNs approximating Kolmogorov PDEs has been given in \cite{deryck2021error}.
Alternative approaches are the deep Ritz method \cite{MR3767958} or the more specific deep Nitsche method which includes essential boundary conditions \cite{MR4247207}.
A comparison of deep Ritz and deep Galerkin has been performed in \cite{MR4199960}.
As an alternative to these approaches, Gaussian Process regression has been studied e.g. in \cite{MR3845646,MR3880138}.

The substitute models generated in this way can often be evaluated very quickly, but so far there have hardly been any approaches to
certification or to speed up the training phase, for example through error evaluations that can be evaluated quickly and can be used
for adaptive learning algorithms.
The informative value of the resulting surrogate models for unseen data is therefore often unsatisfactory, so that a reliable use of
such surrogate models in critical areas has so far not been possible.

For the first time (to the best of our knowledge), we thus propose an adaptive combination of model and data-based reduction processes in a combined method that is efficient and certified at the same time.
We recently presented first ideas of such a combined approach in
the context of transport-diffusion-reaction systems in a
non-adaptive \cite{GHI+2021} and non-certified \cite{HOS2021} fashion.
Using machine learning approaches such as deep neural networks (DNNs) or
kernel methods (KMs), we will  be able to generate adaptive and certified surrogate models.

The structure of our article is as follows.
In \cref{sec:abstract_setting} we introduce our new abstract framework for a certified adaptive surrogate modeling based on a hierarchy of full order, reduced order and machine learning models. In the subsequent \cref{sec:approximation_methods} the approach is specified for parametrized parabolic PDEs with linear output functionals and specific choices of model reduction and machine learning approaches. Finally, several numerical experiments are given in \cref{sec:experiments} to demonstrate the efficiency of the approach in different application scenarios.
For a visual representation of the proposed adaptive hierarchy, we recommend a peek at~\cref{fig:algorithm_visualization} in~\cref{sec:abstract_setting}.
We also provide the source code used to carry out the experiments and generate the figures\footnote{\texttt{\url{https://github.com/ftschindler/paper-2022-certified-adaptive-RB-ML-ROM-hierarchy}}}.

\section{Abstract framework: an adaptive hierarchy of certified surrogate models}
\label{sec:abstract_setting}

Given a parameter domain $\Params \subseteq \R^p$ for $p \in \N$ and end-time $T > 0$, we seek to efficiently and accurately approximate the evaluation of a quantity of interest (QoI) $$f\coloneqq S \circ A\colon\Params \to L^2([0, T]),$$ which is not given explicitly but by applying an output operator $$S\colon L^2(0, T; V) \to L^2([0, T])$$ to an intermediate state-trajectory obtained from a state-evaluation operator $$A\colon \Params \to L^2(0, T; V),$$ with a (real-valued) Hilbert space $V$
(we formalize this structure in the notion of a \emph{state-based model} in \cref{def:state_based_model} below).
For simplicity, we restrict ourselves to time-invariant non-parametric output operators, consider a fixed temporal approximation space $Q_\Dt([0, T]) \subset L^2([0, T])$ of finite (but possibly large) dimension $K \coloneqq \dim Q_\Dt([0, T])$ (induced by a temporal grid $\{0 \eqqcolon t_1< \dots < t_K \coloneqq T\} \subset [0, T]$) and state-trajectories in Bochner-type spaces, i.e.~$Q_\Dt(0, T; V) \subset L^2(0, T; V)$.
The presented concepts straightforwardly generalize to time-dependent parametric output operators
and to other spatio-temporal approximations under suitable assumptions.
The subscript $*$ below will be specified later to
distinguish different models to approximate $f$.
As a particular instance of such models, we will specify our approach for parametrized parabolic PDEs with linear output functionals in \cref{sec:approximation_methods} below.

\begin{definition}[{State-based $\textsc{Model}[V_*, A_*, S_*]$}]
  \label{def:state_based_model}
  A state-based model, abbreviated by $M_* \coloneqq \textsc{Model}[V_*, A_*, S_*]$, of order $N_*\coloneqq \dim V_*$ is given by a real-valued Hilbert space $V_*$ (with inner product $(\cdot,\cdot)_{V_*}$ and induced norm $\|\cdot\|_{V_*}$) as (approximate) \emph{state-space}, an (approximate) \emph{state-evaluation operator} $A_*\colon \Params \to Q_\Dt(0, T; V_*)$ and an (approximate) \emph{output operator} $S_*\colon Q_\Dt(0, T; V_*) \to Q_\Dt([0, T])$, such that we obtain the (approximate) QoI $f_*\coloneqq S_* \circ A_*\colon \Params \to Q_\Dt(0, T)$.
  For algorithmic use, each state-based model $M_*$ provides the routines
  \begin{itemize}
    \item $u_*(\mu) \gets M_*\code{.eval_state}[\mu]$, yielding a state-trajectory $u_*(\mu) \coloneqq A_*(\mu)$, and
    \item $f_*(\mu) \gets M_*\code{.eval_output}[\mu]$, yielding the models output $f_*(\mu) \coloneqq S_*\big(A_*(\mu)\big)$,
  \end{itemize}
  each given an input $\mu \in \Params$.
\end{definition}

Such models arise naturally, for instance, when considering time-dependent PDEs and combining the method of lines with a spatial discretization
(corresponding examples are given in \cref{sec:approximation_methods}).
Though we consider an adaptive hierarchy of models in this work, we restrict ourselves to models which each have
a fixed state-space.
While adaptive Finite Element methods relying solely on refinement could be cast into the presented framework, we are rather interested in the situation where we are given a fixed reference model (a common assumption in the field of model order reduction), and refer to \cite{OS2015,ASU2016,Yano16,MR3880338}
for the more general case.

\begin{assumption}[The (costly) ``truth'' FOM]
  \label{asmpt:FOM}
  We assume we are given a state-based model $M_h\coloneqq \textsc{Model}[V_h, A_h, S_h]$ with fixed approximate state-space $V_h$, state-evaluation and output operators $A_h$ and $S_h$, respectively, such that the quality of approximate state-trajectories $M_h\code{.eval_state}[\mu]$ as well as approximate quantities of interest $M_h\code{.eval_output}[\mu]$ is sufficient for all $\mu \in \Params_\textnormal{outer-loop}$, where
  $\Params_\textnormal{outer-loop} \subseteq \Params$
  denotes the (possibly unknown) set of parameters which are of interest for the application at hand.
  We call this model the \emph{full order model (FOM)} as it will be the reference for all other subsequent models.
  We also assume $M_h$ to be costly in the sense that it is not possible to evaluate $M_h$ for all $\mu \in \Params_\textnormal{outer-loop}$ within the given time- and/or resource-constraints, while few evaluations of $M_h$ are feasible.
\end{assumption}

In practice, the sequence $\Params_\textnormal{outer-loop} \subseteq \Params$ could for instance be generated by
a random sampling strategy for uncertainty quantification or by
a parameter vector trajectory obtained by iterative optimization methods. Both these application settings will be
addressed in our experiments.

Since we seek \emph{efficient} as well as \emph{accurate} approximations of $M_h$, we consider two distinct classes of surrogates, each with their respective shortcomings, in the following two subsections.  Thereafter, we propose a combined approach to overcome these deficits.

\subsection{Accurate (linear) reduced order models (ROMs)}
\label{sec:abstract_RB}

Inspired by pro\-jec\-tion-based MOR, we consider struc\-ture-preserving ROMs which mimic the state-based model from \cref{def:state_based_model}.
As it turns out, at least in the context of our envisioned model hierarchy, this notion is also appropriate for ML-based surrogates.

\begin{definition}[{Certified reduced order model $\textsc{ROM}[V_*, A_*, S_*, E_*]$}] \ \\
\label{def:certified_ROM}\noindent
A certified reduced order model (ROM), abbreviated by $M_* \coloneqq \textsc{ROM}[V_*, A_*, S_*, E_*]$, to approximate the FOM $M_h$, is given by a reduced state-space $V_* \subset V_h$, a reduced state-evaluation operator $A_*\colon\Params\to Q_\Dt(0, T; V_*)$, and a reduced output operator $S_*\colon Q_\Dt(0, T; V_*) \to Q_\Dt(0, T)$, such that $M_*$ is itself a state-based $\textsc{Model}[V_*, A_*, S_*]$ of order $N_*\coloneqq \dim V_*$ in the sense of \cref{def:state_based_model}, and that $u_*(\mu)\coloneqq A_*(\mu)$ is an approximation of $u_h(\mu)=A_h(\mu)$ and $f_*(\mu)\coloneqq S_*\big(u_*(\mu)\big)$ is an approximation of $f_h(\mu) = S_h\big(u_h(\mu)\big)$.

  In addition, a certified ROM is supposed to come with an output-error estimation operator $E_*\colon Q_\Dt(0, T; V_*) \times \Params \to \R$ quantifying the \emph{output-approximation error}, i.e.\
  \begin{align*}
    \|f_h(\mu) - S_*(v_*)\|_{L^2([0, T])} &\leq E_*(v_*; \mu), &&\text{for }v_* \in V_*, \mu \in \Params.
  \end{align*}
  For algorithmic use, a ROM provides the routine
  \begin{itemize}
    \item $\Delta_*^f(\mu) \gets M_*\code{.est_output}[\mu]$, yielding the output error estimate
    $E_*\big(u_*(\mu); \mu\big)$,
  \end{itemize}
  given an input $\mu \in \Params$, in addition to the routines from \cref{def:state_based_model}.
\end{definition}

Such ROMs have the potential to allow for an efficient and accurate approximation of the QoI $f_h$.
However, a closer look at how these ROMs are constructed is in order, to make the notions of efficiency and accuracy more precise.
Note that the ROMs considered here only reduce the spatial space $V_h$ and do not reduce the temporal complexity.
The abstract framework, however, could also be applied to space-time ROMs (cf. \cite{MR3166968,MR3208882,MR4249844}).

\begin{assumption}[Aspects of ROMs]
  \label{asmpt:ROM_aspects}\! \!\!
  Given a ROM $M_* \coloneqq \textsc{ROM}[V_*, A_*, S_*, E_*]$ to approximate the FOM $M_h$, we identify three core aspects of $M_*$.
  \begin{enumerate}
    \item\label{asmpt:ROM_aspect:subspace}
      The linear subspace $V_* \subset V_h$ defines the approximation properties of $M_*$. The accuracy is thus limited by the
      Kolmogorov width of the solution manifold~\cite{OhlRav16}.
      We require $N_* \ll N_h$, which usually can only be achieved by training in a problem-adapted (often iterative) manner.
    \item\label{asmpt:ROM_aspect:precomputation}
      We assume that it is possible to carry out pre-computations to prepare the reduced state- and output operators $A_*$ and $S_*$ \dots
    \item\label{asmpt:ROM_aspect:evaluation}
      \dots such that the computational complexity of invoking $M_*\code{.eval_state}$,\linebreak $M_*\code{.eval_output}$, and $M_*\code{.est_output}$ depends on $N_*$ and not on $N_h$.
  \end{enumerate}
\end{assumption}

Traditionally, aspects \ref{asmpt:ROM_aspects}.\ref{asmpt:ROM_aspect:subspace} and \ref{asmpt:ROM_aspects}.\ref{asmpt:ROM_aspect:precomputation} are addressed by a computationally expensive offline-phase, where the ROM is built. On the other hand, aspect \ref{asmpt:ROM_aspects}.\ref{asmpt:ROM_aspect:evaluation} allows for a computationally efficient online-phase, where the ROM is repeatedly evaluated.
The shortcomings of such a fixed offline/online-decomposition are the missing flexibility in optimizing the overall computational costs based on the demands of an outer algorithm and
the missing adaptive choice of accuracy. To overcome these shortcomings, adaptive enrichment has been suggested to incrementally refine initially coarse ROMs (cf.   \cite{MR3431132,MR4014784,MR3702345,zbMATH07358914,Zahr2015,La2014}).
Thus, we aim at not employing ROMs directly, but rather to leverage (a hierarchy of) ROMs to provide computationally
efficient algorithms yielding efficient approximations of prescribed accuracy in an adaptive manner.
We therefore take a closer look at prototypical examples of how to obtain ROMs, before formalizing the
notion of a ROM generator in \cref{def:ROM_generator} below.

\begin{example}[Offline training vs.~online-adapted ROM]
  \label{ex:offline_vs_online_adaptive}
  Let $M_h$ denote the FOM from \cref{asmpt:FOM} and let $G_* \coloneqq \textsc{RomGen}_*[M_h, \varepsilon]$ denote an $\varepsilon$-accurate ROM generator in the sense of \cref{def:ROM_generator} below, for a prescribed accuracy $\varepsilon > 0$.
  Let $\Params_\textnormal{outer-loop} \subset \Params$ denote a set of parameters of interest, which we wish to query the ROM for.%
  \vspace{-2ex}\begin{multicols}{2}
    {\small{\normalfont\runinhead{Goal-oriented offline greedy training.}}%
    For a finite set $\Params_\textnormal{t} \subset \Params$, we train offline:}%
    {\normalfont\begin{algorithmic}[1]
      \State $M_* \gets G_*\code{.precompute}[\,]$
      \While{\hS{-1}$\max\limits_{\mu \in \Params_\textnormal{t}}M_*\code{.est_output}[\mu]\hS{-1}> \varepsilon$\hS{-1}}
      \State $\overline{\mu} \gets \underset{\mu \in \Params_\textnormal{t}}{\arg\max}\, M_*\code{.est_output}[\mu]$
        \State $G_*\code{.extend}[\overline{\mu}]$
        \State $M_* \gets G_*\code{.precompute}[\,]$
      \EndWhile
    \end{algorithmic}}%
    {\small\noindent Online, we query the ROM by:}
    {\normalfont\begin{algorithmic}[1]
      \ForAll{$\mu \in \Params_\textnormal{outer-loop}$}
        \State $f_*(\mu) \gets M_*\code{.eval_output}[\mu]$
      \EndFor
    \end{algorithmic}}%
    {\small{\normalfont\runinhead{Online-adapted ROM.}}\hfill We adaptively train and query the ROM at the same time:}
    {\normalfont\begin{algorithmic}[1]
      \State $M_* \gets G_*\code{.precompute}[\,]$
      \ForAll{$\mu \in \Params_\textnormal{outer-loop}$}
        \If{$M_*\code{.est_output}[\mu] > \varepsilon$}
          \State $G_*\code{.extend}[\mu]$
          \State $M_* \gets G_*\code{.precompute}[\,]$
        \EndIf
        \State $f_*(\mu) \gets M_*\code{.eval_output}[\mu]$
      \EndFor
    \end{algorithmic}}%
    \vspace{6ex}\mbox{}
  \end{multicols}
\end{example}

The key ingredient to building ROMs in an iterative manner is thus a \emph{generator}, which is responsible for covering aspects \ref{asmpt:ROM_aspects}.\ref{asmpt:ROM_aspect:subspace} and \ref{asmpt:ROM_aspects}.\ref{asmpt:ROM_aspect:precomputation} by keeping track of the training inputs, by generating the required training data using the reference model, and by precomputing the quantities which are required to approximate
the state-evaluation and output operators.
The selection of the training inputs, however, is not part of the generator and left to either an iterative algorithm as above, or an a priori choice.
In what follows, similarly to $*$, $\diamond$ will be used as a placeholder subscript indicating different ROM types.

\begin{definition}[{$\varepsilon$-accurate ROM generator $\textsc{RomGen}_\diamond[M_*, \varepsilon]$}]
  \label{def:ROM_generator}
  Given a state-based model $M_* = \textsc{Model}[V_*, A_*, S_*]$ and a target accuracy $\varepsilon > 0$, a ROM generator $G_\diamond \coloneqq \textsc{RomGen}_\diamond[M_*, \varepsilon]$ is given by an initially empty set of training parameters $\Params_\diamond = \emptyset$,
  and the routines
  \begin{itemize}
    \item $G_\diamond\code{.extend}[\mu]$, collecting the input $\Params_\diamond \gets \Params_\diamond \cup \{\mu\}$ and
      any additional data associated with $\mu$ (e.g.~by evaluating $M_*$) required for \code{precompute}, given an input $\mu \in \Params$, and
    \item $M_\diamond \gets G_\diamond\code{.precompute}[\,]$, yielding a certified ROM $M_\diamond\coloneqq \textsc{ROM}[V_\diamond, A_\diamond, S_\diamond, E_\diamond]$ based on the
      data collected
      so far that is \emph{$\varepsilon$-accurate} w.r.t.~the training data, as in
      \begin{align}
        M_\diamond\code{.est_output}[\mu] \leq \varepsilon &&\text{for all }\mu\in\Params_\diamond.
        \label{eq:output_reproduction}
      \end{align}
  \end{itemize}
  Optionally, a generator may provide the routine
  \begin{itemize}
    \item $\widetilde{G}_\diamond \gets G_\diamond\code{.prolong}[\widetilde{M_*}]$,
      yielding a new ROM generator $\widetilde{G}_\diamond \coloneqq\textsc{RomGen}_\diamond[\widetilde{M_*}, \varepsilon]$, given an extended state-based model $\widetilde{M_*} = \textsc{Model}[\widetilde{V_*}, \widetilde{A_*}, \widetilde{S_*}]$, by prolonging all previously collected data onto an extended state-space $\widetilde{V_*}\supset V_*$.
  \end{itemize}
\end{definition}

Specific instances of such generators will be given later in \cref{sec:approximation_methods}. For now it is sufficient to
understand that the \code{extend} routine for an RB model could
be the computation of a new trajectory for the new parameter.
The \code{precompute} method could be a subsequent reduced basis extension for an RB model, or
a retraining of an ML-based surrogate.
The \code{prolong} method can be required if the existing training data for an ML-based surrogate
needs to be extended/padded by zeros for retraining,
e.g. if the dimensionality of the target RB-space has increased.

Using the structure of state-based models we may thus formulate algorithms to generate and
query ROMs to approximate the QoI as \emph{accurately} as desired.
While such approaches have the potential to be computationally
\emph{efficient} compared to repeatedly querying the FOM, there is still room for further improvement.
For example, time-marching integration schemes for ROMs typically need to iteratively
perform nontrivial operations (e.g.\ solving linear systems for implicit solvers)
in each time step. This iterative dependence of the computations make these models still
expensive for $K \gg 1$ time-steps.
Furthermore, ROMs suffer inherently from the scarce-data scenario, i.e.\ the FOM data is very expensive and
typically only few reliable FOM solves can be afforded for generating a ROM.

Machine Learning based approximate models are a possible remedy for this shortcoming.

\subsection{Efficient (nonlinear) machine-learning surrogate models}
\label{sec:abstract_ML}

Broadly speaking, supervised machine learning surrogates aim to
learn an approximate map $\varPhi\colon X \to Y$ from given training data
$(x_i,y_i) \in X\times Y$.
While in general the training data is assumed to follow some unknown random distribution,
in the surrogate modeling context, the target data typically is
deterministically dependent on the inputs, i.e.\ $y_i= \widetilde{\varPhi}(x_i)$ for
some unknown function $\widetilde{\varPhi}\colon X \rightarrow Y$.
A straightforward approach would be to directly learn the FOM output
map $f_h\colon \R^p \to \R^{K}$ by querying~$M_h$.
Due to \cref{asmpt:FOM}, however, we have put ourselves in a position where it is computationally
infeasible to query $M_h$ sufficiently often.
In addition, to the best of our knowledge, there is no computable a posteriori estimate on the generalization
error of machine learning surrogates available.
This prevents the notion of a certified machine-learning surrogate model and thus drastically limits the
usefulness of these techniques if certification is required.
To remedy the issue of lacking training data, another approach would
thus be to learn the ROM output map $f_\ROM\colon \R^p \to \R^{K}$ by querying the associated
certified ROM $M_\ROM$ \cite{GHI+2021}, which is computationally feasible.
One could even compute upper bounds on the QoI error by means of the ROMs output estimate and triangle
inequality \cite{HOS2021}.
However, these bounds are overly pessimistic and the computational cost of certification would involve
the evaluation of a ROM, discarding the expected computational advantage of ML-models.
In general, there is also no hope of using the precomputed quantities required to evaluate the
certificate $M_\ROM\code{.est_output}[\mu]$ to build a certificate for the envisioned ML model,
since these usually require an associated ROM state, which a direct approximation of $f_\ROM$ cannot provide.

\subsection{An accurate and efficient adaptive model}
\label{sec:abstract_adaptive_model}

We thus propose to employ machine-learning to build a certified ML-ROM mimicking the structure of state-based models.
Namely, given a certified ROM $M_* = \textsc{ROM}[V_*, A_*, S_*,$ $E_*]$, we suggest to learn an approximation $A_\ML$ of the reduced state-evaluation operator $A_*$ to obtain another certified ROM $M_\ML\coloneqq \textsc{ROM}[V_*, A_\ML, S_*, E_*]$.
While inheriting the approximation properties and certification qualities of the underlying ROM $M_*$, the ML-based ROM $M_\ML$
may overcome the inherent iterative time-dependency of $M_*$ when evaluating the model, which tackles one of the ROMs shortcomings identified above.
As with any other ROM, we assume the training of the ML-based state-evaluation operator~$A_\ML$ to be carried out by means of a ROM-generator in the sense of \cref{def:ROM_generator}, say $\textsc{RomGen}_\ML$, based on data computed using the underlying ROM $M_*$ (see \cref{sec:approximation_methods} for concrete examples).
Combining these different surrogate models, we finally propose an adaptive model which unites the accuracy and certification of RB ROMs with the efficiency of ML-based surrogates.

\begin{proposition}[{$\varepsilon$-accurate \textsc{AdaptiveModel}$[\textsc{RomGen}_\RB, \textsc{RomGen}_\ML, \varepsilon]$}]
  \label{lem:adaptive_output_model} \ \\
  Let $M_h$ denote the FOM from \cref{asmpt:FOM}, let $\varepsilon > 0$ denote a desired accuracy of the QoI-approximation, and let $\textsc{RomGen}_\RB$ and $\textsc{RomGen}_\ML$ denote $\varepsilon$-accurate ROM-generators in the sense of \cref{def:ROM_generator}, the latter supporting \code{prolong}ation.

  Abbreviated by $M_\textnormal{adapt}\coloneqq\textsc{AdaptiveModel}[\textsc{RomGen}_\RB, \textsc{RomGen}_\ML, \varepsilon]$, we consider an \emph{adaptive $\varepsilon$-accurate model} to approximate FOM states and outputs, which adaptively builds ROMs as required,
  (based on initially trivial ROMs and generators $G_\RB\coloneqq \textsc{RomGen}_\RB[M_h, \varepsilon]$, $M_\RB \gets G_\RB\code{.precompute}[\,]$, $G_\ML \coloneqq \textsc{RomGen}_\ML[M_\RB, \varepsilon]$ and $M_\ML \gets G_\ML\code{.precompute}[\,]$),
  by specifying the routines
  \begin{itemize}
    \item $f_\textnormal{adapt}(\mu) \gets M_\textnormal{adapt}\code{.eval_output}[\mu]$, given by \cref{alg:adaptive_output}, and
    \item $u_\textnormal{adapt}(\mu) \gets M_\textnormal{adapt}\code{.eval_state}[\mu]$, analogous to \cref{alg:adaptive_output},
  \end{itemize}
  each given an input $\mu \in \Params$.
  We then have
  \begin{align}
    \|f_h(\mu) - M_\textnormal{adapt}\code{.eval_output}[\mu]\|_{L^2([0, T])} &\leq \varepsilon &&\text{for all } \mu \in \Params.
    \label{eq:accuracy_of_adaptive_model}
  \end{align}
\end{proposition}

\begin{figure}[ht]
	\centering%
	\footnotesize%
	\includegraphics[width=\textwidth]{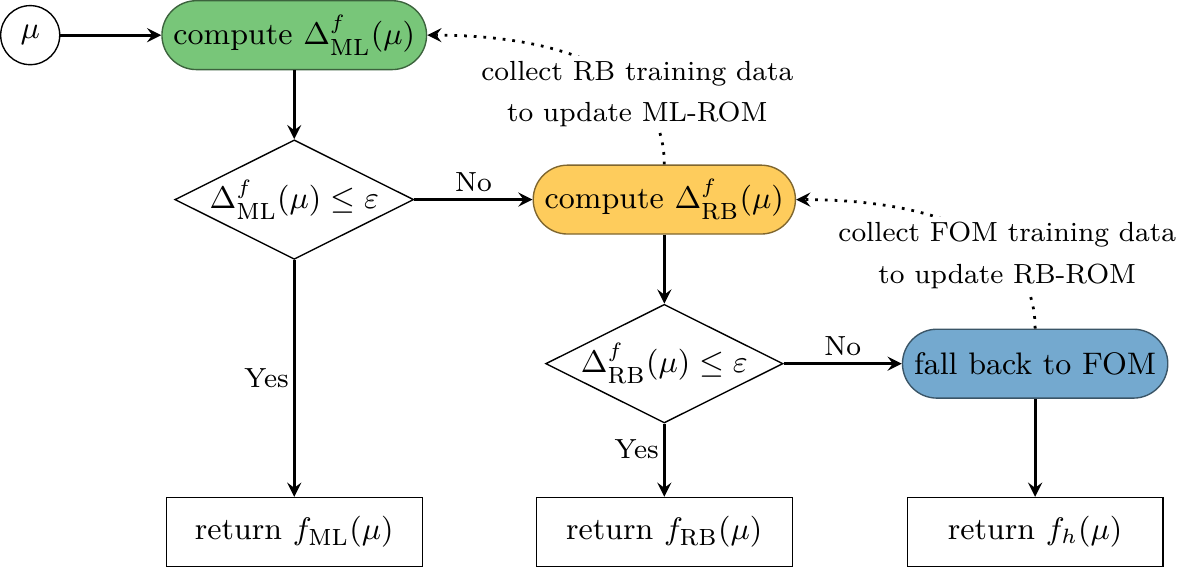}
	\caption{%
		Flow diagram of the adaptive $\varepsilon$-accurate \code{eval_output} from~\cref{alg:adaptive_output}.
	}\label{fig:algorithm_visualization}
\end{figure}

The above model is constructed to produce $\varepsilon$-accurate QoI-approximations with a minimum amount of FOM-evaluations, see~\cref{fig:algorithm_visualization} for a control flow diagram.
While the guaranteed accuracy \eqref{eq:accuracy_of_adaptive_model} follows from \cref{alg:adaptive_output}, we demonstrate its desired performance in actual experiments in \cref{sec:experiments}.
Note that, while we envision RB- and ML-based ROMs, the proposed adaptive model actually holds for any combination of surrogates which fit into this sections' abstract framework.

\begin{algorithm}[h]
  \caption{Adaptive $\varepsilon$-accurate \code{eval_output} from \cref{lem:adaptive_output_model}.}
  \label{alg:adaptive_output}
  With the assumptions and notation from \cref{lem:adaptive_output_model}, the following algorithm adaptively enriches the ROMs $M_\RB$ and $M_\ML$ employing their respective generators $G_\RB$ and $G_\ML$, to provide an output $f_\textnormal{adapt}(\mu)$ of prescribed accuracy.
  \begin{algorithmic}[1]
    \State $\Delta_\ML^f(\mu) \gets M_\ML\code{.est_output}[\mu]$
    \If{$\Delta_\ML^f(\mu) \leq \varepsilon$}\Comment{The ML-prediction is good enough.}
      \State $f_\ML(\mu) \gets M_\ML\code{.eval_output}[\mu]$
      \State \Return $f_\textnormal{adapt}(\mu)\coloneqq f_\ML(\mu)$
    \Else\Comment{The ML-prediction was not good enough, fall back to the RB model.}
      \State $\Delta_\RB^f(\mu) \gets M_\RB\code{.est_output}[\mu]$
      \If{$\Delta_\RB^f(\mu) \leq \varepsilon$}\Comment{The RB-prediction is good enough.}
        \State $G_\ML\code{.extend}[\mu]$\Comment{Collect RB data for ML training.}
        \Statex (optionally:) $M_\ML \gets G_\ML\code{.precompute}[\,]$\Comment{Fit ML with updated data.}
        \State $f_\RB(\mu) \gets M_\RB\code{.eval_output}[\mu]$
        \State \Return $f_\textnormal{adapt}(\mu)\coloneqq f_\RB(\mu)$
      \Else\Comment{Neither ML nor RB prediction was good enough, fall back to the FOM.}
        \State $G_\RB\code{.extend}[\mu]$\Comment{Collect FOM data for RB training.}
        \State $M_\RB \gets G_\RB\code{.precompute}[\,]$
        \State $G_\ML \gets G_\ML\code{.prolong}[M_\RB]$\Comment{Update collected ML training data.}
        \State $G_\ML\code{.extend}[\mu]$\Comment{Collect RB data for ML training.}
        \State $M_\ML \gets G_\ML\code{.precompute}[\,]$\Comment{Fit the ML model anew.}
        \State $f_\RB(\mu) \gets M_\RB\code{.eval_output}[\mu]$
        \State \Return $f_\textnormal{adapt}(\mu)\coloneqq f_\RB(\mu)$
      \EndIf
    \EndIf
  \end{algorithmic}
\end{algorithm}

\section{Certified surrogates for parametrized parabolic PDEs}
\label{sec:approximation_methods}

The methodology presented in the previous section is applicable to any state-based model, for which a reduction scheme with a posteriori error estimates is available.
Here, we focus on models arising from applying the method of lines to linear parabolic PDEs.
Posed on a bounded Lipschitz-domain $\Omega \subset \R^d$, with a given Gelfand triple of suitable Hilbert spaces (encoding essential boundary conditions) $V \subset H^1(\Omega) \subset L^2(\Omega) \subset V'$, we consider state trajectories $u(\mu) \in L^2(0, T; V)$ with $\partial_t u(\mu) \in L^2(0, T; V')$, such that
\begin{align}
  \left<\partial_t u(\mu), v\right> + a\big(u(\mu), v; \mu\big) &= l(v; \mu) &&\text{for all } v \in V,\label{eq:weak_formulation}\\
  u(0; \mu) &= u_0(\mu), &&\text{for } u_0(\mu) \in V,\notag
\end{align}
as well as derived outputs $f\colon \Params \to L^2(0, T)$, $f(t; \mu) \coloneqq s\big(u(t; \mu)\big)$.
Here, for each $\mu \in \Params$, we denote by $l(\,\cdot\,; \mu), s \in V'$ continuous linear functionals and by $a(\,\cdot\,,\,\cdot\,; \mu)\colon V \times V \to \R$ a continuous and coercive bilinear form with symmetric and non-symmetric parts $a_\text{s}, a_\text{n}\colon V \times V \to \R$, such that $a(u, v; \mu) = a_\text{s}(u, v; \mu) + a_\text{n}(u, v; \mu)$ for all $u, v \in V$.
As inner product, we choose $(u, v)_V \coloneqq a_\text{s}(u, v; \overline{\mu})$ for some fixed $\overline{\mu} \in \Params$, and presume a parametrization
of \eqref{eq:weak_formulation}, such that $a$ is continuous and uniformly coercive w.r.t~$(\cdot, \cdot)_V$ in the sense that there exists $\underline{\alpha} > 0$ and
a strictly positive $\alpha\colon \Params \to \R$, $\alpha \geq \underline{\alpha}$, such that $\alpha(\mu)\,\|u\|_V^2 \leq a(u, u; \mu)$ for all $u \in V$
and all $\mu \in \Params$.

\begin{example}[The method of lines]
  \label{ex:fom_method_of_lines}
  To approximate \eqref{eq:weak_formulation}, we use a state-based model in the sense of \cref{def:state_based_model} by employing the method of lines based on a conforming spatial approximation (for simplicity), induced by a Finite Element space $V_h \subset V$ of finite dimension $N_h\coloneqq\dim V_h$, inheriting the product and norm of $V$.
  To obtain the FOM $M_h\coloneqq\textsc{Model}[V_h, A_h, S_h]$ from \cref{asmpt:FOM}, we assume $V_h$ to be suffiently rich, and specify
  \begin{itemize}
    \item
      as \emph{state-evaluation operator} $A_h\colon \Params \to Q_\Dt(0, T; V_h)$, $\mu \mapsto A_h(\mu) \coloneqq u_h(\mu)$, by specifying $u_h(\mu; t_1) = \Pi_{V_h}\big(u_0(\mu)\big) \in V_h$ and $u_h(\mu; t_k) \in V_h$ for all $1 < k \leq K$ as the solution of
      \begin{align}
        (m + \Delta t\, a)\big(u_h(\mu; t_k), v; \mu\big) &= m(u_h(\mu; t_{k - 1}), v; \mu\big) + \Delta t\, l(v; \mu),
      \label{eq:parabolic_FOM_solution}
      \end{align}
      for all $v \in V_h$, where $m(u, v; \mu) \coloneqq (u, v)_{L^2}$ and $\Pi_{V_h}\colon V \to V_h$ denotes a suitable projection operator; and
    \item
      as \emph{output operator} $S_h\colon Q_\Dt(0, T; V_h) \to Q_\Dt(0, T)$, by specifying\\ $S_h(v_h; t_k) \coloneqq s_h\big(v_h(t_k)\big)$ for $1 \leq k \leq K$ and $v_h \in Q_\Dt(0, T; V_h)$;
  \end{itemize}
  such that we obtain the FOM QoI $f_h\colon \Params \to Q_\Dt(0, T)$ as $f_h(\mu) \coloneqq \big(S_h \circ A_h\big)(\mu)$.
\end{example}

\subsection{Tier 1: linear subspaces and certification by RB methods}
\label{sec:RB}

As special instantiation of a certified reduced order model, we will
make use of a Reduced Basis (RB) ROM for parabolic problems.
RB-methods are meanwhile standard for model reduction of a variety of types
of parametric FOMs \cite{Hesthaven2016,Haasdonk2017} as they both can provide
accurate approximation with low-dimensional ROMs, as well as certification by
residual-based a posteriori error estimators.
These methods comprise a so called offline-phase, where
a reduced space is determined by invoking the FOM for several
parameter realizations, gathering so called snapshots and compressing these
to a reduced basis spanning the RB space $V_\RB$.
This basis generation is typically done in a greedy fashion, i.e.\
iteratively determining a new parameter instance for which the
current model is providing the worst approximation,
invoking the FOM for this parameter and extending the reduced basis
based on this FOM state information.
Frequently, a compression of (error) trajectory data is performed by a
so called proper orthogonal decomposition (POD) \cite{MR1868765}
or an approximate variant called hierarchical approximate
POD (HAPOD) \cite{MR3860899}.
For time-dependent problems, this then results in the POD-Greedy procedure,
which both works very well in practice, and has provable
excellent convergence rates \cite{HO08,H13}.
After computing the reduced basis, the offline-phase allows to compute
further quantities such as projected system components for the ROM.
After single execution of this offline-phase, in the subsequent online-phase
the reduced model can be assembled and solved in an online-efficient
fashion for arbitrary parameter realizations.
Online-efficiency here refers to a computational complexity which is
independent of the spatial dimension of the FOM.

In view of \cref{sec:abstract_setting}, employing RB-methods has two components: a certified ROM in the sense of \cref{def:certified_ROM}, the purpose of which is to perform the online-efficient computations, and a ROM generator in the sense of \cref{def:ROM_generator}, the purpose of which is to generate the reduced space and to precompute all quantities required for an efficient evaluation of the ROM.
We specify both further below, based on a Galerkin-projection of the FOM onto the reduced space.
However, we first detail a central aspect of RB models, namely certification by error bounds.

In the context of the FOM from \cref{ex:fom_method_of_lines}, suppose we are given a reduced solution $u_\RB(\mu) \in V_\RB \subset V_h$ of \eqref{eq:parabolic-RB-weak-form} as an approximation of the solution $u_h(\mu) \in V_h$, for $\mu \in \Params$.
Inserted in the FOM \eqref{eq:parabolic_FOM_solution}, it will deviate from equality
by a residual term~$R_k$. Interestingly, the norm of these residuals can be computed even without knowing
the FOM solution.
For parabolic problems, we employ the a posteriori error estimate from \cite{GP05,Haasdonk2017} to obtain the bound
\begin{align}
  \|u_h(\mu) - u_\RB(\mu)\|_{L^2(0, T; V_h)} &\leq E_\RB^u\big(u_\RB(\mu); \mu\big)
\label{eq:rb_state_estimate}\\
                                             &\coloneqq \frac{1}{\alpha(\mu)} \Big(\sum_{k=1}^{K-1}\Delta t\|R_k\big(u_\RB(t_k; \mu); \mu\big)\|_{V_h'}^2\Big)^{1/2},
\notag
\end{align}
if $u_h(0; \mu) \in V_\RB$.
Owing to our choice of the energy norm, $\alpha$ can be efficiently and explicitly computed by means of the min-$\theta$ approach \cite{Haasdonk2017} for a large class of affinely decomposed problems, or estimated by the more involved successive constraint minimization method \cite{HUYNH2007} in general.
Overall, all quantities determining $E_\RB^u$ can be computed online-efficient
in complexity independent of the dimension of $V_h$, cf.~\cite{Haasdonk2017} and also \cite{BEOR2014} for a numerically stable projection for the computation of the dual norm of the residual.
Using the linearity and continuity of the output functional, which gives the bound
\begin{align}
  |f_h(t_k; \mu) - f_\RB(t_k; \mu)| &\leq |s\big(u_h(t_k; \mu) - u_\RB(t_k; \mu); \mu\big)|
\notag\\
                                    &\leq \|s(\cdot; \mu)\|_{V_h'}\; \|u_h(t_k; \mu) - u_\RB(t_k; \mu)\|_{V_h}
\notag
\end{align}
  point-wise in time, together with the same arguments as for \eqref{eq:rb_state_estimate}, we obtain the following statement.

\begin{proposition}[Efficiently computable RB-output error estimate \cite{GP05}]
  \label{lem:abstract_rb_output_estimate}
  Let~$f_h$ denote the FOM output from \cref{ex:fom_method_of_lines} and let $E_\RB^u$ denote the a posteriori RB-state estimate from \eqref{eq:rb_state_estimate}.
  We then have for any $v_\RB\in Q_\Dt(0, T; V_\RB)$
  \begin{align}
    \|f_h(\mu) - S_\RB(v_\RB)\|_{L^2(0, T)} &\leq E_\RB\big(v_\RB; \mu\big) \coloneqq \|s(\cdot; \mu)\|_{V_h'}\, E_\RB^u\big(v_\RB; \mu\big),
    \label{eq:abstract_rb_output_estimate}
  \end{align}
  where $S_\RB$ denotes the reduced output operator from \cref{ex:RB_ROM}.
\end{proposition}

We can thus estimate the output-error induced by any element of the reduced space $V_\RB$, which motivates the design of the ML-ROM further below.
In particular, we readily obtain an estimate for the output $f_\RB(\mu)$ of the RB-ROM from \cref{ex:RB_ROM}, induced by the corresponding reduced solution trajectory $u_\RB(\mu$):
\begin{align}
  \|f_h(\mu) - f_\RB(\mu)\|_{L^2(0, T)} &\leq E_\RB\big(u_\RB(\mu); \mu\big) &&\text{for all } \mu \in \Params.
\label{eq:rb_output_estimate}
\end{align}
These estimates will be the basis for the $\varepsilon$-accurate RB-ROM generator from \cref{ex:POD_RB_generator} and the following ROM.

\begin{example}[RB-Galerkin-ROM for the method of lines]
  \label{ex:RB_ROM}
  Given the FOM $M_h = \textsc{Model}[V_h, A_h, S_h]$ from \cref{ex:fom_method_of_lines}, and given a reduced basis space $V_\RB \subset V_h$ of dimension $N_\RB\coloneqq\dim V_\RB$, by a suitable generator, we obtain the RB-ROM for the method of lines, $M_\RB\coloneqq\textsc{ROM}[V_\RB, A_\RB, S_\RB, E_\RB]$, as a certified ROM  in the sense of \cref{def:certified_ROM} by specifying as output-estimation operator the one from \eqref{eq:rb_output_estimate}, and by specifying
  \begin{itemize}
    \item
      as reduced \emph{state evaluation operator} $A_\RB\colon \Params \to Q_\Dt(0, T; V_\RB)$ with $\mu \mapsto A_\RB(\mu)\coloneqq u_\RB(\mu)$,
      the iterative Galerkin-projection of \eqref{eq:parabolic_FOM_solution} onto $V_\RB$, i.e.\
      $u_\RB(\mu; t_1) \coloneqq\linebreak \Pi_{V_\RB} u_h(\mu; t_1)$ and $u_\RB(\mu; t_k) \in V_\RB$ for $1 < k \leq K$ as the solution of
      \begin{align}
        \label{eq:parabolic-RB-weak-form}
        (m + \Delta t\, a)\big(u_{\RB}(\mu; t_k),v;\mu\big) = m\big(u_{\RB}(\mu; t_{k - 1}),v;\mu\big) + \Delta t\, l(v;\mu)
      \end{align}
      for all $v\in V_{\RB}$, where $\Pi_{V_\RB}\colon V_h \to V_\RB$ denotes a suitable projection (e.g.\ $L^2$-orthogonal projection);
    \item
      as reduced \emph{output operator} $S_\RB\colon Q_\Dt(0, T; V_\RB) \to Q_\Dt(0, T)$, by specifying\\ $S_\RB(v_\RB; t_k) \coloneqq s_h\big(v_\RB(t_k); \mu\big)$ for $1 \leq k \leq K$ and $v_\RB \in Q_\Dt(0, T; V_\RB)$;
  \end{itemize}
  such that we obtain the approximate QoI $f_\RB\colon \Params \to Q_\Dt(0, T)$ as $f_\RB(\mu) \coloneqq \big(S_\RB \circ A_\RB\big)(\mu)$.
\end{example}

The online-efficiency of the above model stems from the fact that the restriction of the FOM state-evaluation and output operators to the reduced subspace can be precomputed (see the generator below), such that the computational complexity of their evaluation in the above model depends on $N_\RB$ instead of $N_h$.
We proceed similarly to the basis extension of the POD-Greedy algorithm \cite{H13}, while keeping more information from the added trajectory in view of the output-reproduction \eqref{eq:output_reproduction}.
While any POD-algorithm could be used in the generator, in our context of high temporal resolution $K\gg 1$, we employ the hierarchical approximate POD (HaPOD) \cite{MR3860899} to allow for data compression already during data generation, on the fly.

\begin{example}[HaPOD RB-Galerkin-ROM generator]
  \label{ex:POD_RB_generator}
  Given the FOM $M_h = \textsc{Model}$ $[V_h, A_h, S_h]$ from \cref{ex:fom_method_of_lines} and a target-accuracy $\varepsilon > 0$, we obtain an $\varepsilon$-accurate ROM generator in the sense of \cref{def:ROM_generator} to generate the RB-Galerkin-ROM from \cref{ex:RB_ROM}, by hierarchically building a reduced basis $\varphi_\RB^{(0)} \subset \varphi_\RB^{(1)}\subset\dots$ with $\varphi_\RB^{(0)}\coloneqq\emptyset$ by means of the HaPOD \cite{MR3860899}.
  Starting with an empty set of collected training parameters $\Params_\RB^{(0)} \coloneqq\emptyset$, we define:
  \begin{itemize}
    \item
      As a means to \code{extend} the basis for a new parameter $\mu \in \Params$, we set $n\coloneqq |\Params_\RB|$,
      collect $\Params_\RB^{(n+1)} = \Params_\RB^{(n)} \cup \{ \mu \}$, and extend the reduced basis with the remains of the trajectory $u_h(\mu)\coloneqq M_h\code{.eval_state}[\mu]$, after projection onto the current reduced space $V_\RB^{(n)}\coloneqq \langle\varphi_\RB^{(n)}\rangle$, compression by the \code{HaPOD}-algorithm and re-orthonormaliza\-tion w.r.t.~$(\cdot,\cdot)_{V_h}$ by a given \code{gram_schmidt}-algorithm,
      \begin{align}
        \varphi_\RB^{(n+1)} \coloneqq \code{gram_schmidt(}\varphi_\RB^{(n)} \cup \code{HaPOD(} u_h(\mu) - \Pi_{V_\RB^{(n)}}u_h(\mu), \varepsilon_\textnormal{pod} \code{))},
      \notag
      \end{align}
      where $\Pi_{V_\RB^{(n)}}$ denotes a $(\cdot,\cdot)_{V_h}$-orthonormal projection onto $V_\RB^{(n)}$ and where the POD-tolerance $\varepsilon_\textnormal{pod} > 0$ needs to be chosen such that the estimate from \cref{lem:abstract_rb_output_estimate} applied to the solution $u_\RB^{(n + 1)}(\mu) \in V_\RB^{(n + 1)}$ of the RB-ROM obtained using the extended reduced space $V_\RB^{(n+1)}\coloneqq \langle\varphi_\RB^{(n+1)}\rangle$ indicates output-reproduction in the sense of \eqref{eq:output_reproduction} up to the required accuracy, i.e.~$E_\RB\big(u_\RB^{(n + 1)}(\mu); \mu\big) < \varepsilon$ (e.g.~by setting $\varepsilon_\textnormal{pod}$ to machine accuracy).
    \item
      As a means to \code{precompute} all data required for an online-efficient RB-ROM, we employ the standard precomputation by a Galerkin projection of the parameter-separable components of $m$, $a$ and $l$ from \eqref{eq:parabolic_FOM_solution} onto the current reduced space $V_\RB\coloneqq \langle\varphi^{(n)}\rangle$ and refer to \cite[Cor.~2.29]{Haasdonk2017} and the monographs mentioned in the introduction for details and generalizations.
      Given the precomputed data, we obtain the online-efficient RB-ROM from \cref{ex:RB_ROM} that can be evaluated and certified with a computational complexity depending only on the dimension of the current reduced space $V_\RB$, in contrast to the dimension of $V_h$.
  \end{itemize}
  Since the reduced basis is built hierarchically and since the output-reproduction is ensured in each step, the resulting RB-ROM is $\varepsilon$-accurate in the sense of \cref{def:ROM_generator}.
\end{example}

Note that by using incremental POD algorithms such as the HaPOD \cite{MR3860899}, the above generator may \code{extend} the given reduced basis with a FOM trajectory $u_h(\mu)$, without actually storing the full trajectory.

\subsection{Tier 2: non-linear efficient evaluation by ML methods}
\label{sec:certified_approaches}

As discussed in \cref{sec:abstract_setting}, the evaluation of the RB-ROM $M_\RB$ suffers from the iterative nature of the time-dependency of the reduced state-evaluation operator $A_\RB$, which is bound to the same time-stepping scheme as the FOM one.
While not an issue for small problems, this bottle-neck limits the use of RB-ROMs in many-query contexts where the underlying problem requires long-time integration or a prescribed temporal accuracy, yielding a large number of time-steps $K \gg 1$.
As motivated in \cref{sec:abstract_adaptive_model}, we thus seek to obtain a certified ML-based ROM by replacing the reduced state-evaluation operator $A_\RB$ with a non-linear ML-based prediction.
While many use-cases for ML-based predictions
are conceivable, we restrict ourselves to two obvious choices which are
also evaluated numerically in \cref{sec:experiments}.
We postpone the definition of the corresponding generators to the following subsections, where we also detail various concrete ML realizations.
The two following models simply presume an already learned approximate map $\varPhi$, as introduced in \cref{sec:abstract_ML}, and can be used to speed up any state-based ROM, not just the one from \cref{ex:RB_ROM}.

\begin{example}[Certified ML-ROMs for instationary problems]\label{ex:certified-ml-roms}\ \\
  Let $M_* = \textsc{ROM}[V_*, A_*,$ $S_*, E_*]$ denote a certified ROM of order $N_*\coloneqq\dim V_*$.
  We obtain a certified ML-ROM, $M_\ML = \textsc{ROM}[V_*,$ $A_\ML, S_*, E_*]$, by specifying the ML-reduced state-evaluation operator $A_\ML\colon\Params \to Q_\Dt(0, T; V_*)$, $\mu \mapsto u_\ML(\mu)$ by one of the following means:

  \emph{The time-``vectorized'' case:} let $\varPhi\colon\R^p \to \R^{K\times N_*}$ be a ML-based map, trained by a suitable generator, to predict the ROM coefficients for all time-steps at once.
  We then obtain the ML-state-evaluation $u_\ML(\mu)$ by
  \begin{align}
    \underline{u_\ML(\mu; t_k)}_n = {\varPhi(\mu)}_{k, n}, &&\text{for }1 \leq n \leq N_*\text{ and }1 \leq k \leq K,
    \label{eq:vectorized_ML_model}
  \end{align}
  where $\underline{u_\ML(\mu; t_k)} \in \R^{N_*}$ denotes the DoF-vector of $u_\ML(\mu; t_k) \in V_*$.

  \emph{The ``random-access''-in-time case:} let $\varPhi\colon\R^{p + 1} \to \R^{N_*}$ be a ML-based map, trained by a suitable generator, to predict the ROM coefficients at any point in time.
  We then obtain the ML-state-evaluation $u_\ML(\mu)$ by
  \begin{align}
    \underline{u_\ML(\mu; t_k)}_n = {\varPhi(\mu, t_k)}_n, &&\text{for }1 \leq n \leq N_*\text{ and }1 \leq k \leq K.
    \label{eq:random_access_ML_model}
  \end{align}
\end{example}

Note that, while the evaluation of the reduced state-evaluation operator $A_*$ has been replaced by the ML-reduced operator $A_\ML$ in the above example, the cost of evaluating the output $S_*$ as well as the a posterior error estimator $E_*$ remains unchanged.
However, both of the latter only involve forward-applications of functionals (corresponding to a cost of $\mathcal{O}(K\cdot N_*)$, which can also be carried out partly in parallel w.r.t. $K$), while the former involves an incremental inversion of $K$ operators (corresponding to a cost of $\mathcal{O}(K\cdot N_*^2)$).
Thus, the computational gain in actual experiments may be significant (as demonstrated below).

\subsubsection{HaPOD-VKOGA-ROM}
\label{sec:hapod_vkoga_generator}

In the following we briefly introduce Kernel methods (KMs), which revolve around the notion of (strictly) positive definite kernels.
They are a mathematically very well analyzed set of techniques and algorithms and are widely applied in the context of numerical mathematics, approximation theory and machine learning such as surrogate modeling \cite{Santin2021}.
For more details, we refer the interested reader to the standard textbooks \cite{Wendland2005, Fasshauer2007, Fasshauer2015}.

For a given set $\Omega_\text{d}$ of input data, where we assume $\Omega_\text{d} \subset \R^d$ for simplicity,  a scalar-valued, strictly positive definite kernel is a symmetric function $k_s\colon \Omega_\text{d} \times \Omega_\text{d} \rightarrow \R$, such that the so called kernel matrix $(k_s(x_i, x_j))_{i,j=1}^n$ is positive definite for any choice of pairwise distinct data points $x_i\in\Omega_\text{d}$, $i=1,\dots, n$, with $n \in \N$ arbitrary.
As a well known example of a strictly positive definite kernel, we employ the Gaussian kernel $k_s(x, y) = \exp(-\Vert x - y \Vert_2^2)$ in our experiments in \cref{sec:experiments}.
Kernel methods can be used for unsupervised as well as supervised learning tasks.
In the case of supervised regression, machine learning usually considers input data
$X_N\coloneqq \{x_1,\dots, x_N\} \subset \Omega_\text{d}$ and corresponding target values $\{y_1,\dots, y_N\} \subset \R$.
In the case of a mean-square error (MSE) loss, one usually seeks for
\begin{align}
  \varPhi = \mathrm{argmin}_{f \in \mathcal{H}}\ \mathcal{L}(f),&&\text{with}&&\mathcal{L}(f) = \frac{1}{N} \sum_{i=1}^N |y_i - f(x_i)|^2 + \lambda \cdot \Vert f \Vert_{\mathcal{H}}^2, \label{eq:ERM}
\end{align}
where $\lambda \geq 0$ is a regularization parameter to avoid overfitting to possibly noisy target data.
If kernels are used for this task, the space $\mathcal{H}$ is chosen as a Reproducing Kernel Hilbert Space w.r.t.~$k_s$,
and a kernel representer theorem states that the solution of \eqref{eq:ERM} can be written as
\begin{align}
\label{eq:kernel_model}
\varPhi^{(N)}(\cdot) = \sum_{i=1}^N \alpha_i^{(N)} k_s(\cdot, x_i),
\end{align}
i.e.\ in a data dependent way.
The coefficients $\alpha_i \in \R$ can be obtained by solving an $N \times N$ linear equation system, which is however only efficient for small or moderate data sizes $N$.
Especially, the resulting model can be computed exactly and is unique, which is a major advantage over gradient based optimization algorithms
which are used for neural networks.
Further advantages of kernel methods are their appealing mathematical theory including error estimates, optimality and uniqueness statements, which allow to analyze the approximation properties of the kernel models. \\
In order to obtain a sparse kernel model one can choose a suitable subset
$X_n \subset X_N$ of $n \ll N$ training data points.
For a computationally tractable selection one usually employs greedy selection strategies, which are well analyzed in terms of their convergence rates \cite{SH16b, wenzel2021analysis}.
For this, one usually makes use of a change of basis, which allows to efficiently update the coefficients of the kernel model and thus alleviates the need for a recomputation, which would be necessary when using the standard representation of \eqref{eq:kernel_model}.
This efficient stepwise update can also be used in the case of newly incoming data, as required in the experiments in \cref{sec:num_exp_optimization}. \\
To also learn vector valued data in $\R^r$, one uses matrix valued kernels, which are an extension of the concept of scalar valued kernels.
For our experiments, we make use of matrix-valued kernels of the form $k\colon \Omega_\text{d} \times \Omega_\text{d} \rightarrow \R^{r \times r}$, $k = k_s \cdot I_r$, which is a suitable choice if no prior knowledge about symmetries or constraints is available.
In this case, the kernel model from \eqref{eq:kernel_model} has vector valued coefficients $\alpha_i^{(N)} \in \R^r$.
We denote those greedy algorithms using vector-valued kernels by VKOGA \cite{Santin2021}.

To obtain a VKOGA-based generator in the sense of \cref{def:ROM_generator} to be used in \cref{alg:adaptive_output}, we specify its \code{extend}, \code{precompute}, and \code{prolong} methods:
\begin{itemize}
\item \code{extend}: The \code{extend} method adds samples of the form
  $\big(\mu, \underline{u_\RB(\mu;t_k)}_{k=1, ..., K} \big) \in \Params \times\R^{N_* \cdot K}$
to the training set, where $\underline{u_\RB(\mu;t_k)}_{k=1, ..., K}$ is obtained via the\\\code{eval_state} method of the RB-based ROM.
\item \code{precompute}: The \code{precompute} method iteratively updates the kernel model by
  using the efficient greedy implementation of VKOGA to obtain $\varPhi$.
The resulting kernel model $\varPhi$ is then used for the ML-state-evaluation, in particular in the vector-valued
time-vectorized prediction from \cref{ex:certified-ml-roms} (which
gives both the coefficients and their time evolution at once without the need for an iterative time integration).
\item \code{prolong}:
        The \code{prolong} method is called in the case that the reduced basis had to be extended due to an inaccurate RB-ROM.
	This basis extension necessitates
        an adjustment of the VKOGA model,
        because the output dimension has to be adapted.
	Accordingly, the previously collected training data also has
        to be padded by zeros in order
        to be able to reuse that data.
\end{itemize}

\subsubsection{HaPOD-DNN-ROM}
\label{sec:hapod_dnn_generator}

As an example for the ``random-access''-in-time ML-ROM from \cref{ex:certified-ml-roms}, we consider deep neural networks (DNNs) as surrogates where the time-component is incorporated in the input of the network.
\par
DNNs are machine learning algorithms that recently also gained
attention in the context of reduced order modeling, see for instance~\cite{HU2018,LC2020,WWHZ2021,Fresca2022,WHR2019,KPRS2019}.
In this paper we focus on one particular example of DNNs, namely \emph{feedforward neural networks}, which apply an alternating sequence of affine-linear mappings and element-wise non-linearities to the input \cite{PV2018,EGJS2018}.

	Let therefore $L\in\N$, $N_0,\dots,N_L\in\N$, $W_i\in\R^{N_i\times N_{i-1}}$
	for $i=1,\dots,L$, and $b_i\in\R^{N_i}$ for $i=1,\dots,L$. We call
	$L$ the number of layers in the neural network, and
	$N_0,\dots,N_L$ the number of neurons in each of the layers.
	Further, the entries of the matrices $W_1,\dots,W_L$ are called
	weights, and the entries of the vectors $b_1,\dots,b_L$ are
	called biases. The corresponding neural network is the
	$L$-tuple $\big((W_1,b_1),\dots,(W_L,b_L)\big)$. In addition, let
	$\rho\colon\R\to\R$ be the so called activation function. We
	denote the function defined by the neural network
	$\big((W_1,b_1),\dots,(W_L,b_L)\big)$, using~$\rho$ as activation
	function, by~$\varPhi\colon\R^{N_0}\to\R^{N_L}$. The function~$\varPhi$
	is, for a given input $x\in\R^{N_0}$, defined by
		$\varPhi(x) \coloneqq r^L(x)$,
	where $r^L(x)$ is derived recursively as
		$r^0(x) \coloneqq x$,
		$r^i(x) \coloneqq \rho^*\left(W_i\, r^{i-1}(x) + b_i\right), i=1,\dots,L-1$, and
		$r^L(x) \coloneqq W_L\,r^{L-1}(x) + b_L$.
	Here, the function $\rho^*\colon\bigcup_{d\in\N}\R^d\to\bigcup_{d\in\N}\R^d$ denotes
	the element-wise application of the activation function $\rho$,
	i.e.\ for all $d\in\N$ and $z\in\R^d$ we define
	$\rho^*(z)\coloneqq\left(\rho(z_1),\dots,\rho(z_d)\right)\in\R^d$.

The weights and biases of a neural network are subject to training
by minimizing a loss function that measures the difference between the output produced by the neural network and the target outputs on the training samples
analogous to \cref{eq:ERM}.
Usual choices for the optimization routine are variants of (stochastic) gradient descent, see for instance
\cite{BCN2018} for a survey of optimization algorithms employed in
machine learning.
The gradient of the loss function with respect to the weights and biases of the neural network is efficiently computed using backpropagation~\cite{RHW1986}.
The iterations of the optimizer are called training epochs.
Since the training results typically depend strongly on the initialization of weights and biases of the DNN, multiple training restarts are performed using different initial guesses.
A method to prevent the neural network from overfitting the given training data is to perform early stopping (see for instance~\cite{P1997} for a discussion of different approaches) of the training if the loss on a validation set (usually chosen to be disjoint from the training set) does not decrease for a prescribed number of consecutive training epochs.
\par
Since DNNs are fitted in a supervised manner, they are particularly well-suited for the task of predicting reduced coefficients after being trained on a limited amount of samples.
This idea has first been developed in~\cite{HU2018} for the stationary case and was extended in~\cite{WHR2019} to the instationary case, based on learning the (orthogonal) projection of the FOM trajectories onto the RB space.
Instead, we propose to learn the reduced coefficients from RB trajectories which, in combination with the HaPOD generator from \cref{ex:POD_RB_generator}, eliminates the need for storing or recomputing the FOM trajectory.
\par
To obtain a DNN-based generator in the sense of \cref{def:ROM_generator} to be used in \cref{alg:adaptive_output}, we specify its \code{extend}, \code{precompute}, and \code{prolong} methods:
\begin{itemize}
	\item \code{extend}:
	Within this method, pairs of the form $\big((\mu,t_k),\underline{u_\RB(\mu;t_k)}\big)\in(\Params\times[0,T])\times\R^{N_*}$ for $n=1,\dots,N_*$ are added to the set of training samples.
	Hence, to obtain $\underline{u_\RB(\mu;t_k)}$, the \code{eval_state} routine of the RB-based ROM is called.
	The result is stored and reused in the subsequent call to \code{eval_output} in \cref{alg:adaptive_output}.
	In distinction from the VKOGA-ROM presented above, the number of training samples for the DNN-ROM scales with the number $K$ of time steps, and is thus much larger than for the VKOGA-ROM.
	This constitutes one of the main differences between the time-``vectorized'' case and the ``random-access''-in-time case.
	\item \code{precompute}: In the \code{precompute} method, the neural network is trained as described above.
  The resulting function $\varPhi\colon\R^{p+1}\to\R^{N_*}$ is used for the ML-state-evaluation $A_\ML$ in \cref{ex:certified-ml-roms}.
	To reduce the computational effort, it is beneficial not to train for each new training sample, but to collect a certain amount of data before performing another neural network training.
	This further ensures that a reasonable number of samples is available for training.
	It might also be beneficial to reuse previously trained networks by using their weights and biases as initial guess for the optimizer (this is certainly only possible if the architecture did not change, see the \code{prolong} method).
      \item \code{prolong}: As in the case of VKOGA-ROM, this routine is called after the RB space has been extended. Then the target values of the existing training data
        are padded with zeros in the new dimensions and the neural network architecture is adjusted to match the now extended output-dimension.
\end{itemize}

We should recall at this point that the function $\varPhi$ can be evaluated in a
random-access-in-time manner, in particular
for multiple time instances at the same time.
This constitutes a major difference to the corresponding state evaluation of the RB-based ROMs, where the computation of reduced states
is only possible in a strict iterative and computationally non-ignorable manner due to the dependence of the solution on the previous time step.
The parallelization of computing $\varPhi$ for different time instances thus indicates a major contribution for possible speed up compared to the RB-ROM.

\section{Numerical experiments}
\label{sec:experiments}

As a central point for our implementation, we use the freely available Python-based MOR library pyMOR\footnote{\url{https://pymor.org/}} \cite{MRS2016}:
\begin{itemize}
  \item
    For the FOM from \cref{asmpt:FOM}, we employ pyMORs builtin numpy/scipy-based \code{discretize_instationary_cg} method in \cref{subsec:num_exp_mc} and a similar discretizer based on Python-bindings of the freely available C++ PDE solver library dune-gdt\footnote{\url{https://docs.dune-gdt.org/}} in \cref{sec:num_exp_optimization}, each to obtain pyMORs \code{InstationaryModel}s which match the state-based models from \cref{def:state_based_model}.
  \item
    We employ pyMORs \code{inc_vectorarray_hapod} implementation of the HaPOD, as well as the numerically stable \code{gram_schmidt} to build the RB-ROM generator from \cref{ex:POD_RB_generator}, where we use pyMORs \code{ParabolicRBReductor} to perform the precomputation of the reduced state-evaluation and output operators, as well as the assembly of the a posteriori error estimator.
  \item
    The DNN-based model from \cref{ex:certified-ml-roms} and the associated generator from \cref{sec:hapod_dnn_generator} are realized as pyMOR models and reductors, based on pytorch-lightning\footnote{\url{https://github.com/PyTorchLightning/pytorch-lightning}} and PyTorch \cite{PGM+2019}.
  \item
    The VKOGA-based model from \cref{ex:certified-ml-roms} and the associated generator from \cref{sec:hapod_vkoga_generator} are also realized as pyMOR models and reductors, based on the VKOGA software library \cite{Santin2021}\footnote{\url{https://gitlab.mathematik.uni-stuttgart.de/pub/ians-anm/vkoga}}.
\end{itemize}

Subsequently, we demonstrate the performance of the adaptive model from \cref{lem:adaptive_output_model} in various scenarios.

\subsection{PDE-constrained minimization in reactive flow}
\label{sec:num_exp_optimization}

We consider single-phase reactive flow through an artificial pipe $\Omega\coloneqq [0, 5]\times [0, 1] = \Omega_\textnormal{w} \cup \Omega_\textnormal{c}$, partitioned into a washcoat $\Omega_\textnormal{w}\coloneqq [0, 5]\times [0, h_\textnormal{w}]$ of height $h_\textnormal{w}\coloneqq 0.34$ and a channel $\Omega_\textnormal{c}\coloneqq \Omega\backslash\Omega_\textnormal{w}$, with inflow and outflow boundaries $\Gamma_\textnormal{in}\coloneqq\{x = 0\}\times [h_\textnormal{w}, 1]$ and $\Gamma_\textnormal{out}\coloneqq \{x = 5\}\times [h_\textnormal{w}, 1]$, respectively.
For each parameter $\mu\coloneqq(\textnormal{Da}, \textnormal{Pe})\in\Params$ within the diffusion-dominated regime $\Params\coloneqq [0.01, 10]\times [9, 11]$, we seek the solution $u(\mu) \in L^2(0, T; V_g)$ with
$\partial_t u(\mu) \in L^2(0, T; V_0')$ of the linear advection-diffusion-reaction equation for $T = 5$,
\begin{align}
  \partial_t u(\mu) - \nabla\cdot\big(\kappa\nabla u(\mu)\big) + \textnormal{Pe}\,\nabla\cdot\big(\vec{v}\,u(\mu)\big) + \textnormal{Da}\,\chi_{\Omega_\textnormal{w}}u(\mu) &= 0 &&\textnormal{in }\Omega \times (0, T),
\label{eq:reactive_flow}\\
  \big(\kappa\nabla u(\mu)\big)\cdot n_{\partial\Omega} &= 0 &&\text{on }\Gamma_\textnormal{out},
\notag\\
    u(t=0; \mu) &= \chi_{\Gamma_\textnormal{in}} &&\text{in }\Omega,
\notag
\end{align}
in a weak sense, where $V_g\coloneqq\{v \in H^1(\Omega) \,|\, v|_{\partial\Omega\backslash\Gamma_\textnormal{out}} = g\text{ in the sense of traces}\}$
and the Dirichlet data $g\coloneqq \chi_{\Gamma_\textnormal{in}}$ encodes essential boundary
conditions, $\chi_*$ denote indicator functions of the sets $*$, the advective field is given by $\vec{v}\coloneqq(\chi_{\Omega_\textnormal{c}}, 0)^\top$, and
the diffusion $\kappa$ is given by $1$ in the channel $\Omega_\textnormal{c}$ and by the first component of the permeability field of the Spe10 Model1 benchmark\footnote{\url{https://www.spe.org/web/csp/datasets/set01.htm}}, rescaled to $[0.001, 1]$, in the washcoat $\Omega_\textnormal{w}$.

Picking an extension $\tilde{g} \in H^1(\Omega)$, such that $\tilde{g}|_{\partial\Omega\backslash\Gamma_\textnormal{out}} = g$ in the sense of traces, we shift \eqref{eq:reactive_flow} from the affine subspace $V_g$ to the linear subspace $V\coloneqq V_0$ (which is more suitable for ROMs) and obtain the bilinear form $a$ and linear functional $l$ in \eqref{eq:weak_formulation} as $a(u, v; \mu)\coloneqq\int_\Omega\kappa(\nabla u\cdot\nabla v)\dx + \mu_1\int_\Omega\nabla\cdot(\vec{v}u)v\dx + \mu_2\int_\Omega\chi_{\Omega_\textnormal{w}}u\,v\dx$ and $l(v; \mu)\coloneqq -a(\tilde{g}, v; \mu)$, respectively.
Thus, \eqref{eq:weak_formulation} with $u_0(\mu)\coloneqq 0$ describes the evolution of the shifted trajectory $\overcirc{u}(\mu) \in L^2(0, T; V)$ and we recover the solution of \eqref{eq:reactive_flow} as $u(\mu) = \overcirc{u}(\mu) + \tilde{g}$.
The output functional $s(v)\coloneqq|\Gamma_\textnormal{out}|^{-1}\int_{\Gamma_\textnormal{out}}v\ds$ measures the average concentration at the outflow, inducing as QoI the \emph{break-through curve} $f(t; \mu) = s\big(u(t; \mu)\big)$.

To obtain the FOM from \cref{ex:fom_method_of_lines}, we partition $\Omega$ into equidistant squares and approximate $V$ by bi-linear continuous Finite Elements $V_h$, yielding $N_h = \dim V_h$ spatial degrees of freedom, and employ $K - 1$ implicit Euler steps in \eqref{eq:parabolic_FOM_solution}, where $N_h$ and $K$ are given in the following subsections.
Note that, while artificially constructed, this example features the most relevant aspects of highly resolved single-phase reactive flow \cite{GHI+2021}.

In the context of \eqref{eq:reactive_flow}, we consider the problem of PDE-constrained minimization with box-constraints, where we seek to minimize the $L^\infty$-misfit
\begin{align}
  \min_{\mu \in \Params} \mathcal{J}(\mu)\coloneqq \|\hat{f}_h - f_\textnormal{adapt}(\mu)\|_{L^\infty([0, T])},
\label{eq:minimization_reactive_flow}
\end{align}
given a desired QoI $\hat{f}_h\coloneqq f_h(\hat{\mu})$ with $\hat{\mu}\coloneqq (5.005, 10)$ being the center of $\Params$ (the form of $\hat{f}_h$ is not known by the minimizer).
While we have no theoretical justification, it is visible from \cref{fig:flow_minimization_VKOGA_adaptive_eps} (top left) that the above minimization is feasible, in that a global minimum exists within $\Params$, and neither local minima nor saddle-points can be observed.
We thus use the Nelder-Mead simplex algorithm \cite{GH2010} from \code{scipy.optimize.minimize}\footnote{\url{https://docs.scipy.org/doc/scipy/reference/optimize.minimize-neldermead.html}} to solve \eqref{eq:minimization_reactive_flow}. 

Note that this choice of optimization method is only depicted here as a demonstrator. Gradient or even Hessian based optimization methods are in general more favorable for such classes of PDE constrained optimization problems and have been investigated in the context of a pure FOM - ROM hierarchy with adaptive enrichment in \cite{Zahr2015,MR3716566,MR4269464,BKMOSV22}.
Derivative information for such methods can be either computed with the help of a dual problem or by directly computing sensitivities with respect to the parameters. In both cases a posteriori error estimates are available for the respective ROM approximations and can thus also be used in an extended setting with a full FOM-ROM-ML hierarchy. We also refer to \cite{KKLOO22}, where a pure FOM-ML hierarchy has been used to accelerate PDE-constrained optimization. In this framework, however, as no intermediate ROM was considered, the certification was achieved by employing the FOM in a post-processing step.

For $f_\textnormal{adapt}(\mu)\coloneqq M_\textnormal{adapt}\code{.eval_output}[\mu]$ in \eqref{eq:minimization_reactive_flow}, we employ the $\varepsilon$-accurate adaptive model $M_\textnormal{adapt} = \textsc{AdaptiveModel}[\textsc{RomGen}_\RB, \textsc{RomGen}_\ML, \varepsilon]$ from \cref{lem:adaptive_output_model}, using the HaPOD RB-Galerkin-ROM generator from \cref{ex:POD_RB_generator} with $\varepsilon_\textnormal{POD} = 10^{-15}$ (for output-reproduction) as $\textsc{RomGen}_\RB$, and the time-vectorized VKOGA ML-ROM generator from \cref{sec:hapod_vkoga_generator} as $\textsc{RomGen}_\ML$, a combination which we abbreviate by HaPOD-VKOGA-ROM.
We use a re-training of the ML-ROM in \cref{alg:adaptive_output}, Line 8, after every call to $G_\ML\code{.extend}[\mu]$, which implies a re-fitting of the underlying VKOGA kernel model in each iteration.
While unusual in the context of ML, we expect the overall optimization to be finished after a few hundred iterations already (rendering extensive data-collection before fitting useless) and also experienced exceptional performance of the VKOGA fitting (for this application) compared to, for instance the fitting of a DNN.
Note that, being able to use ML techniques based on very little data is actually one of the main benefits of the certified ML-ROMs we propose, since the prediction quality can always be efficiently certified a posteriori.

The use of the adaptive model in the context of \eqref{eq:minimization_reactive_flow} raises two questions:
\begin{enumerate}
  \item
    What accuracy $\varepsilon > 0$ of the adaptive model is required for the minimization in \eqref{eq:minimization_reactive_flow} to succeed?
    While the estimate from \cref{lem:abstract_rb_output_estimate} could be used to provide $L^2$-bounds, it is unclear a priori how to relate that to \eqref{eq:minimization_reactive_flow} to obtain a rich enough misfit-landscape for the minimizer to succeed.
  \item
    Considering all evaluations of $M_\textnormal{adapt}$ required during the minimization in \eqref{eq:minimization_reactive_flow}, how often are the individual models $M_h$, $M_\RB$ and $M_\ML$ used internally?
    I.e., to what degree are the adaptively enriched surrogate models sufficient for the application at hand, and how does this compare to a minimization employing only $M_h$, or only $M_h$ and $M_\RB$?
\end{enumerate}

\subsubsection{Minimization with an a priori fixed tolerance $\varepsilon$}
\label{sec:num_exp_optimization:fixed_eps}

To answer the first question above, we consider a coarse FOM $M_h$, where the minimization in \eqref{eq:minimization_reactive_flow} can be carried out using a tolerance of $\varepsilon=0$ (which corresponds to using $f_h(\mu)\coloneqq M_h\code{.eval_output}[\mu]$ instead of $f_\textnormal{adapt}$) to obtain a reference.
We then minimize \eqref{eq:minimization_reactive_flow} using the adaptive model $M_\textnormal{adapt}$ for various a priori choices of the tolerance $\varepsilon$ and compare to the reference run, to study the impact of $\varepsilon$ on the convergence of the minimization, as well as the performance of the adaptive model (which we measure for this experiment by counting the evaluations of the submodels used internally in the adaptive model, compare \cref{alg:adaptive_output}).

\begin{table}[ht]
  \centering%
  \footnotesize%
  \begin{tabular}{l|rcss}
    model used in \eqref{eq:minimization_reactive_flow} & \#evals & conv. & \mcol{c}{rel.~min.~err.} & \mcol{c}{rel.~obj.~err.} \\ \hline
    $M_h(\varepsilon=0)$ & 99 & yes & \sci{2.82}{-6} & \sci{2.44}{-8} \\
    \grayout{$M_\textnormal{adapt}(\varepsilon=1.00\hS{-2}\cdot\hS{-2}10^{-1})$} & \grayout{403} & \grayout{no} & \scigray{2.73}{-1} & \scigray{1.85}{-3} \\
    \grayout{$M_\textnormal{adapt}(\varepsilon=1.50\hS{-2}\cdot\hS{-2}10^{-2})$} & \grayout{403} & \grayout{no} & \scigray{2.16}{-1} & \scigray{1.88}{-3} \\
             $M_\textnormal{adapt}(\varepsilon=1.25\hS{-2}\cdot\hS{-2}10^{-2})$ & 121 & yes & \sci{2.53}{-3} & \sci{2.21}{-5} \\
             $M_\textnormal{adapt}(\varepsilon=1.00\hS{-2}\cdot\hS{-2}10^{-2})$ & 126 & yes & \sci{4.45}{-4} & \sci{4.35}{-6} \\
             $M_\textnormal{adapt}(\varepsilon=1.00\hS{-2}\cdot\hS{-2}10^{-3})$ & 97 & yes & \sci{4.03}{-5} & \sci{4.77}{-7} \\
             $M_\textnormal{adapt}(\varepsilon=1.00\hS{-2}\cdot\hS{-2}10^{-4})$ & 105 & yes & \sci{1.13}{-4} & \sci{7.01}{-7} \\
             $M_\textnormal{adapt}(\varepsilon=1.00\hS{-2}\cdot\hS{-2}10^{-5})$ & 96 & yes & \sci{1.98}{-6} & \sci{3.74}{-8}
  \end{tabular}
  \caption{%
    Performance and accuracy of the models detailed in \cref{fig:flow_minimization_VKOGA_fixed_eps} for the experiment from \cref{sec:num_exp_optimization:fixed_eps}: number of model evaluations (\#evals) required to reach convergence (conv.) or divergence (greyed out) as reported by the minimizer, the relative error in the minimizer $\|\hat{\mu} - \mu\|_{l_2}/\|\hat{\mu}\|_{l_2}$ for the final parameter $\mu$ (rel.~min.~err.), and normalized objective value $\mathcal{J}(\mu)/\mathcal{J}(\mu^{(0)})$ (rel.~obj.~err), where $\mu^{(0)}$ denotes the initial guess and $\hat{\mu}$ denotes the optimal parameter (not known by the minimizer).
  }
  \label{tab:flow_minimization_VKOGA_fixed_eps}
\end{table}
\begin{figure}[ht]
  \centering%
  \footnotesize%
  \includegraphics{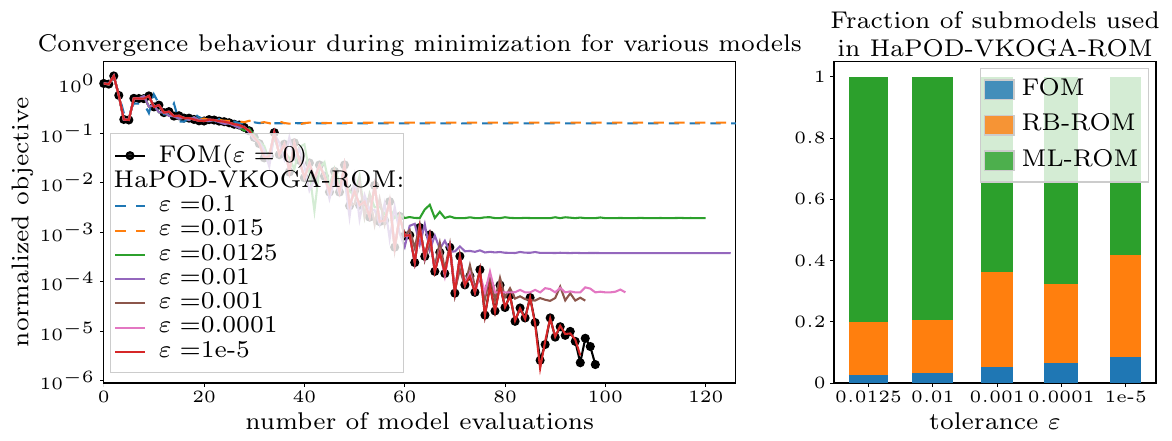}
  \caption{%
    Performance of the adaptive HaPOD-VKOGA-ROM in terms of convergence (left) and use of surrogate models (right) for selected tolerances $\varepsilon$ compared to the FOM, for the experiment from \cref{sec:num_exp_optimization:fixed_eps}, each starting from the same initial guess $\mu^{(0)}\coloneqq (2, 10.5)$, with a coarse FOM ($N_h = \num[group-separator={,}]{2121}$ DoFs, $K = \num[group-separator={,}]{1001}$ time steps).
    Left: evolution of the normalized objective $\mathcal{J}(\cdot)/\mathcal{J}(\mu^{(0)})$ during the minimization in \eqref{eq:minimization_reactive_flow} using the reference FOM $M_h$ (black) or the adaptive model $M_\textnormal{adapt}$ for various tolerances $\varepsilon$ (colored).
    (Note that we depict all intermediate model evaluations during the Nelder-Mead minimization, the $x$-axis thus corresponds to the number of model evaluations, and not the iteration count of the minimization.)
    Note also that, for the adaptive model, convergence of the minimization can only be observed for small enough tolerances (as indicated by solid lines).
    Right: fraction of (surrogate) models used internally in $M_\textnormal{adapt}$ to compute the output $f_\textnormal{adapt}$ in \eqref{eq:minimization_reactive_flow} (compare \cref{alg:adaptive_output}).
  }\label{fig:flow_minimization_VKOGA_fixed_eps}
\end{figure}

As we observe in \cref{tab:flow_minimization_VKOGA_fixed_eps} and \cref{fig:flow_minimization_VKOGA_fixed_eps} (left), the choice of the tolerance $\varepsilon$ has a direct impact on the convergence behaviour of the minimization (if chosen too large, the minimization does not converge), as well as on the quality of the minimum.
On the other hand, the choice of the tolerance also has a direct impact on the performance of the adaptive model: if chosen too small, the computationally cheaper surrogate models are used less often (compare \cref{alg:adaptive_output}), and the adaptive model needs to resort to computationally expensive FOM evaluations (as we observe in \cref{fig:flow_minimization_VKOGA_fixed_eps}, right).
Since neither influence is known a priori, the tolerance $\varepsilon$ has to be chosen by trial-and-error.
Thus, a satisfying answer to the first question from above cannot be given by an a priori choice of $\varepsilon$.

\subsubsection{Minimization with an adaptive tolerance}
\label{sec:num_exp_optimization:adaptive_eps}

To remedy this situation, we propose an adaptive choice of the tolerance during the minimization in the following manner.
Starting from an arbitrarily large initial tolerance $\varepsilon^{(0)} >0$, we compute a running average (of width $n_\textnormal{av.}$) of the objective $\mathcal{J}$, and use the slope of a linear regression of this average over a window of same size to obtain an approximate gradient of $\mathcal{J}$.
Monitoring this approximate gradient, as well as the gradient normalized by the current value of $\mathcal{J}(\cdot)/\mathcal{J}(\mu^{(0)})$, we detect stagnation if the former lies below $\varepsilon_{\nabla\mathcal{J}} > 0$ or the latter lies below $\tilde{\varepsilon}_{\nabla\mathcal{J}} > 0$ for more than $n_\textnormal{stag.}$ consecutive model evaluations.
Upon stagnation, we lower the tolerance of the adaptive model, $\varepsilon^{(n + 1)}\coloneqq \varepsilon^{(n)}/10$, and discard those collected ML training data points, the associated prediction of which does not satisfy the new tolerance.

We thus shifted the issue of choosing an a priori tolerance $\varepsilon$ of the adaptive model, to the specification of the stagnation detection parameters $n_\textnormal{av.}$, $n_\textnormal{stag.}$, $\varepsilon_{\nabla\mathcal{J}}$ and $\tilde{\varepsilon}_{\nabla\mathcal{J}}$.
These can be interpreted with respect to the minimization problem: since the stagnation detection is not applied after each minimization iteration, but after each model evaluation, the averaging window size should be larger than the input dimension (at least for simplex minimizers), e.g.~$n_\textnormal{av.}\coloneqq 2\cdot\dim\Params$, while the role of $n_\textnormal{stag.}$, $\varepsilon_{\nabla\mathcal{J}}$ and $\tilde{\varepsilon}_{\nabla\mathcal{J}}$ is to detect stagnation of the objective before the minimization algorithm falsely reports convergence due to an inaccurate model (e.g.~they can be chosen depending on the tolerances of the minimizer).

\begin{figure}[ht]
  \centering%
  \footnotesize%
  \includegraphics{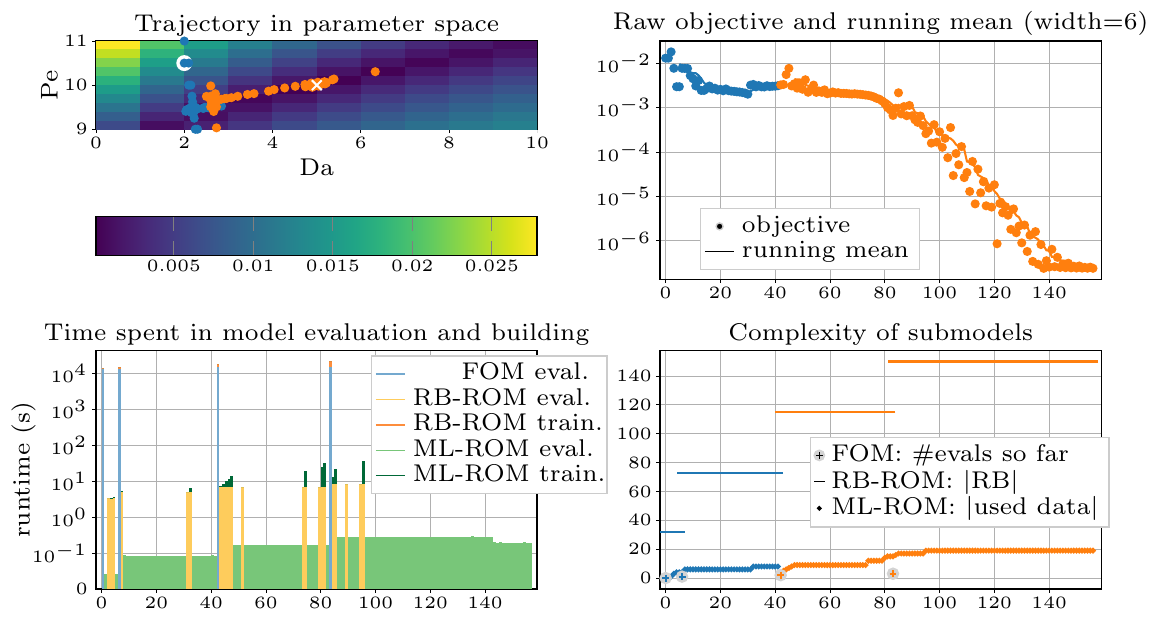}
  \caption{%
    Performance of the adaptive HaPOD-VKOGA-ROM in terms of convergence in parameter space (top left) and objective (top right), and runtime (bottom left) for the experiment from \cref{sec:num_exp_optimization:adaptive_eps} with an adaptive choice of the tolerance $\varepsilon$ (indicated by {\color[rgb]{0.122,0.467,0.706}blue, $\varepsilon^{(0)} = 6.77\cdot 10^{-2}$}, and {\color[rgb]{1.,0.498,0.056}orange, $\varepsilon^{(1)} = 6.77\cdot 10^{-3}$}, in all but the bottom left plot).
    Top left: parameter space $\Params$ and the objective $\mathcal{J}$ obtained from a very coarse model as colored background, initial guess (white circle), found minimum (white cross) and points in parameter space chosen by the Nelder-Mead minimization of \eqref{eq:minimization_reactive_flow}, where the adaptive model $M_\text{adapt}$ was evaluated (colored dots).
    Top right: value of the objective $\mathcal{J}$ for each of the 157 evaluations of the adaptive model $M_\text{adapt}$ required until convergence.
    Bottom left: required time to compute the objective $\mathcal{J}(\mu)$ for each $\mu$ chosen during minimization using $M_\text{adapt}$, and details on time spent for evaluation and building of the required submodels of $M_\text{adapt}$, compare \cref{alg:adaptive_output}.
    Bottom right: complexity of the adaptively generated submodels of $M_\text{adapt}$ in terms of $\dim V_\RB$ for the HaPOD-ROM and number of centers for the VKOGA-ROM (note the drop of the latter due to tolerance adaption).
  }\label{fig:flow_minimization_VKOGA_adaptive_eps}
\end{figure}

We demonstrate the performance of the adaptive model with the proposed adaptive choice of the tolerance $\varepsilon$ for a FOM with $N_h = \num[group-separator={,}]{2051841}$ DoFs and $K = \num[group-separator={,}]{10001}$ time steps, where the computation of $f_h(\hat{\mu})$ required\footnote{On a dual socket compute server equipped with two Intel Xeon E5-2698 v3 CPUs with 16 cores running at 2.30GHz each and 256GB of memory available. Note that the high spatial and temporal resolution of this FOM was picked on purpose: while not required for such an academic example, we intend to mimic more realistic scenarios as motivated in \cite{GHI+2021}.} $1.23\cdot10^4$s (the results are displayed in \cref{fig:flow_minimization_VKOGA_adaptive_eps}).
Regarding the use of the submodels in \cref{alg:adaptive_output}, only four evaluations of $M_h$ are required, 22 evaluations of $M_\RB$ and 131 evaluations of $M_\ML$.
The size of the reduced basis (compare \cref{ex:POD_RB_generator}) grew from 0 to 32, 73, 115 and 150 after the four evaluations of $M_h$.
Starting from an initial tolerance corresponding to $\|\hat{f}\|_{L^2([0, T])}$, the stagnation detection with $n_\textnormal{av.} = 6$, $n_\textnormal{stag.} = 10$, $\varepsilon_{\nabla\mathcal{J}} = -1\cdot 10^{-15}$ and $\tilde{\varepsilon}_{\nabla\mathcal{J}} = 5\cdot 10^{-5}$ correctly identified stagnation due to lacking accuracy after 39 evaluations, a single adaptation of the tolerance was sufficient for convergence.
Overall, the minimization required 157 evaluations of the adaptive model and completed in $23$h and $56$m (corresponding to seven FOM evaluations).
As visible from \cref{fig:flow_minimization_VKOGA_adaptive_eps} (bottom left), while the evaluation of $M_\RB$ is roughly three orders of magnitudes faster than the evaluation of $M_h$, using the ML model is again one to two orders of magnitudes faster than the RB model.
Note that since the minimization terminates rather quickly, the computational gain of employing the ML model is of the same order as the effort required to train it.
In hindsight, restricting the hierarchy to FOM and RB model would have thus been enough, which is, however, not known a priori.
In total, the use of the adaptive model hierarchy renders the minimization of \eqref{eq:minimization_reactive_flow} feasible within a day, while 157 FOM evaluations would have taken some 22 days.

\subsection{Monte Carlo} \label{subsec:num_exp_mc}
As a second application for our proposed algorithm, we consider the
estimation of an average and variance of an output quantity over the parameter domain with
respect to a given probability density function.
To this end, we make use of a Monte Carlo approach that typically requires a
large amount of output evaluations for different parameters.
\par
As before, let $\Params \subseteq \R^p$ denote the parameter domain.
Further, assume that we are given a probability density function
$\rho\colon\Params\to[0,\infty)$ with respect to the standard Lebesgue
measure on $\R^p$, i.e.\ it holds $\int_\Params\rho(\mu)\dmu = 1$. We
are interested in the expected value and the variance of the average of the output
$f(\mu)$ over a time interval $\tau \subseteq (0,T)$ with respect to the density
$\rho$ on the parameter space. Consequently, we seek for
\begin{align*}
	\Expec{\bar{f}} \coloneqq \int_\Params \rho(\mu)\bar{f}(\mu)\dmu\quad\text{and}\quad\Var{\bar{f}} \coloneqq \int_\Params \rho(\mu)\left(\bar{f}(\mu)-\Expec{\bar{f}}\right)^2\dmu,
\end{align*}
where $\bar{f}\colon\Params\to\R$ denotes for a given parameter
$\mu\in\Params$ the average of $f(\mu)$ over time, i.e.~
$\bar{f}(\mu) \coloneqq |\tau|^{-1}\int_\tau f(\mu;t)\dt$.
To estimate the expected value $\Expec{\bar{f}}$ and the variance $\Var{\bar{f}}$, we randomly sample a
set of parameters $\Params_\mc\subseteq\Params$ that is distributed
according to the density $\rho$. The expected value and the variance are then
approximated by unbiased estimators given as
\begin{align*}
	\Expec{\bar{f}} \approx \theta_{\bar{f}} \coloneqq \frac{1}{N_\mc}\sum_{\mu\in\Params_\mc} \bar{f}(\mu)
	\quad\text{and}\quad
	\Var{\bar{f}} \approx \frac{1}{(N_\mc-1)}\sum_{\mu\in\Params_\mc} \left(\bar{f}(\mu)-\theta_{\bar{f}}\right)^2,
\end{align*}
where $N_\mc\coloneqq\card{\Params_\mc}$ denotes the number of drawn
samples.
\par
\begin{figure}
	\begin{minipage}[c]{0.6\textwidth}
    \includegraphics[width=\textwidth]{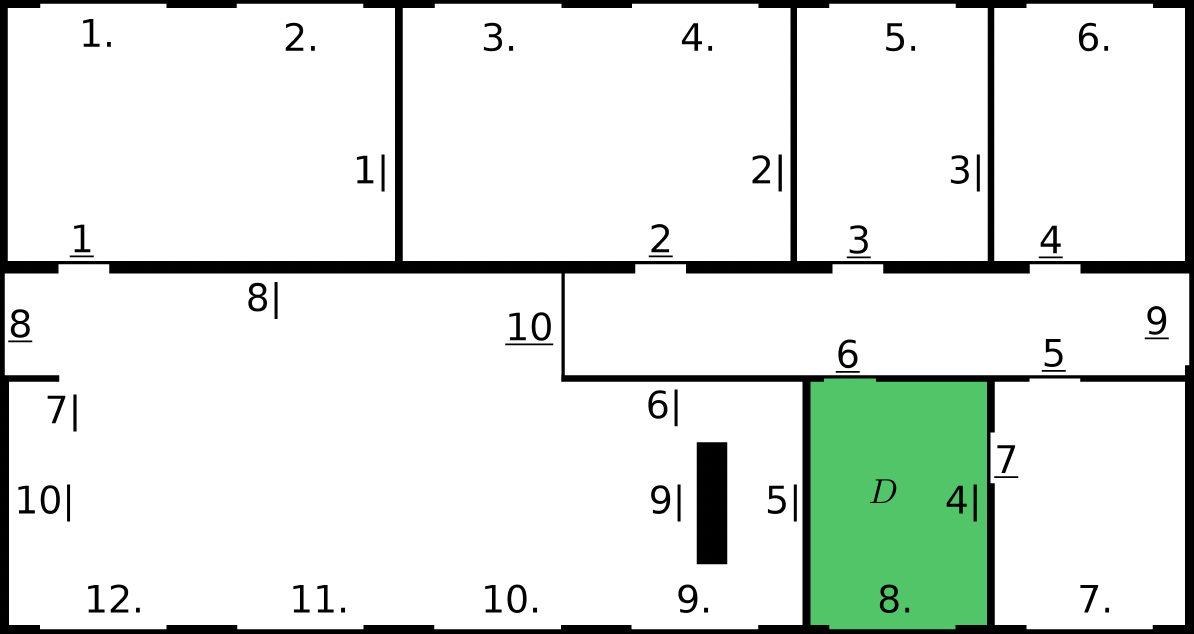}
	\end{minipage}%
	\hfill
	\begin{minipage}[c]{0.375\textwidth}
  \caption{Data functions and domain of interest $D\subset\Omega$ (green) from \cite{MR4269464} for the Monte Carlo experiment from \cref{subsec:num_exp_mc}: building components inducing parametric data functions are denoted as $i|$ (walls), $i.$ (heaters) and $\underline{i}$ (doors).
In this experiment, $7|$, $10|$, $\underline{8}$ and $\underline{9}$ are kept non-parametric.}
		\label{fig:mmexc}
	\end{minipage}
\end{figure}

In our numerical experiment, we consider the average temperature in a
certain room of a building. A similar (stationary) problem in the
context of PDE-constrained optimization has been investigated
in~\cite{MR4269464}. The spatial domain $\Omega \coloneqq (0,2)\times(0,1)\subset\R^2$
describes the building floor, and the room we are interested in is
denoted by $D\subset\Omega$, see also \Cref{fig:mmexc}. The quantity of interest (QoI)
$f$ for a parameter $\mu\in\Params$ at time $t\in[0,T]$ is
given as
$f(\mu; t) \coloneqq \abs{D}^{-1}\int_D u(t,x;\mu) \dx$,
i.e.\ the average temperature at time $t$ over the domain $D$. We
consider the heat equation on the domain $[0,T]\times\Omega$ for $T = 1$ with
parametric diffusion coefficient $\kappa(\mu)\in L^\infty(\Omega)$, $\kappa(\mu) > 0$ for $\mu \in \Params$  almost everywhere,
parametric and time-dependent right hand side source term $\ell(\mu;t)\in L^2(\Omega)$ for $\mu \in \Params$ and $t \in [0, T)$, and
homogeneous Dirichlet boundary conditions on the whole boundary~$\partial\Omega$ of the domain~$\Omega$, and the constant function
$u_0(\mu)\equiv 0$ as initial condition. The bilinear form $a$ and the
linear functional $l$ from \eqref{eq:weak_formulation} are thus given
for $\mu\in\Params$, $t\in[0,T]$, and $u,v\in V$ as
\begin{align*}
	a(u,v;\mu) = \int_\Omega \kappa(\mu)\,\nabla u(x)\cdot\nabla v(x)\dx,\qquad l(v;\mu;t) = \int_\Omega v\,\ell(\mu;t)\dx.
\end{align*}
The $28$-dimensional parameter domain $\Params$ consists of the
diffusion parameters for the $8$ inner walls and $8$ inner doors
(involved in the coefficient function $\kappa$), and the heat emitted
by the $12$ heaters at the top and the bottom of $\Omega$ (part of the
source term $\ell$), cf.\ \Cref{fig:mmexc}.
The inner walls and doors are treated as volumes in the interior of the domain $\Omega$ and are not part of the boundary.
The final source term $\ell$ is multiplied by $\min\{2t,1\}$ to simulate the process of slowly turning on the heaters.
\par
As time interval of interest, we choose $\tau=(0.9,1)$, i.e.\ we consider the last $10\%$ of the whole time domain and measure the average temperature over that time span within the room $D$.

For this example with a higher dimensional input space, we employ neural networks
to learn the required input-output relations as elaborated in
\Cref{sec:certified_approaches}.
The setup of the used neural network is given by a standard feedforward network with 4 hidden layers of size
128 each.
The input dimension is given by $29$ ($28$ parameters and $1$ time variable) and the output dimension by the
dimension of the current RB space. This realizes the ``random-access''-in-time from \cref{ex:certified-ml-roms}.
Since the RB space is enlarged during runtime, also the architecture is
correspondingly adapted by changing the output dimension.
The common ReLU activation function is used for all but the final layer.
This setup is a standard architecture suited for the problem at hand in view of the number of training points that will be available in this Monte Carlo example.
The input
data is standardized to a range of $[-1,1]$ by means of a linear transformation prior to the training.
In contrast, no scaling is applied to the outputs since this would remove the inherent weighting of the individual reduced components.
For the training of the network the Adam optimizer is applied to minimize the MSE-loss.
In addition, we employ an initial learning rate of $5 \cdot 10^{-3}$ and a mini-batch size of~$128$.
The learning rate is modified according to a step learning rate scheduler, which reduces the learning rate by a factor of $0.7$ every $10$ epochs.
\par
The optimization of the ML-ROM was performed for up to $100$ epochs, but early stopping based on the decay of the validation loss (computed on $5\%$ of the overall available data that serves as validation set) within $10$ subsequent epochs was used.
No restarts were performed, because the used network is large enough and the batch size of $128$ is small
in comparison to the amount of training data.
Therefore the stochastic optimization introduces
enough randomness which makes it highly unlikely to get stuck in a bad local minimum.
Furthermore, we are interested in keeping the computational effort for the neural network training low, since in our adaptive algorithm the time required for preparing the reduced order models contributes to the overall computation time.
Due to the certification of the results of the ML-ROM and the possibility of falling back to the RB-ROM, it is moreover not necessary to train a ``perfect'' DNN, but a network that is used sufficiently often to replace the RB-ROM.
\par
For the experiment, $N_\mc=\num[group-separator={,}]{10000}$ randomly sampled parameters were drawn from a uniform distribution on the parameter space $\Params$.
Further, the FOM consists of $N_h=\num[group-separator={,}]{321206}$
DoFs and $K - 1$ time steps were performed with $K=\num[group-separator={,}]{1000}$ . A fixed absolute tolerance of $\varepsilon = 5 \cdot 10^{-2}$ for the adaptive model from \Cref{lem:adaptive_output_model} was used throughout the whole computation.
The training batch size was set to $200$, i.e.\ after $200$ new ROM solutions the optimization of the DNN was restarted (compare line 8 of \cref{alg:adaptive_output}).
\par
\begin{figure}
	\centering
	\includegraphics{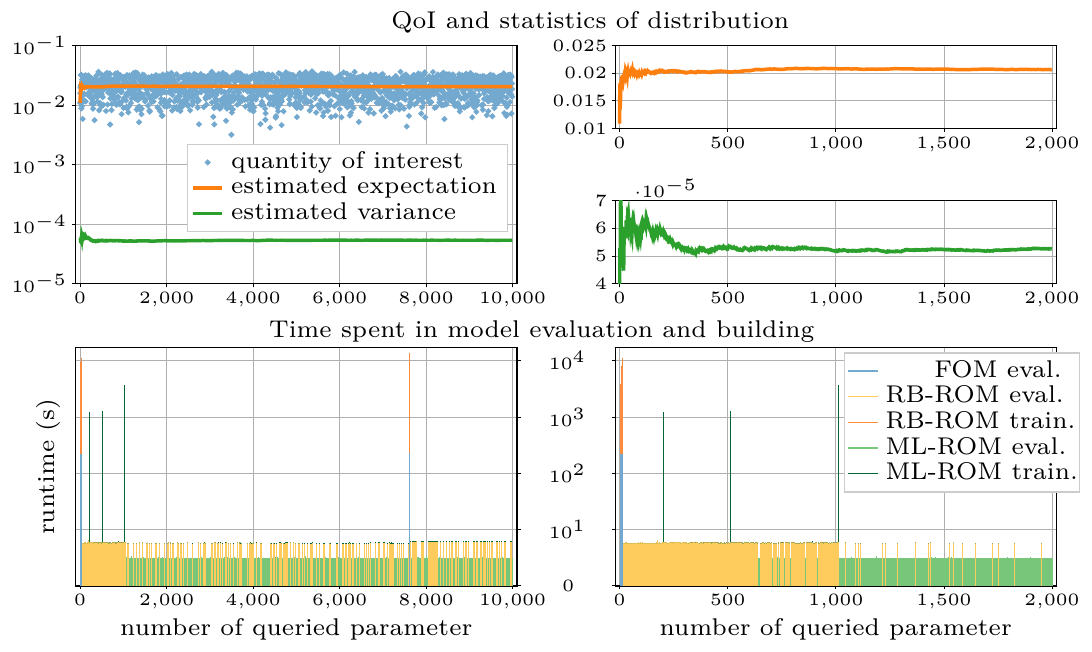}
  \caption{Performance of the adaptive HaPOD-DNN-ROM in terms of Monte Carlo estimation (top row) and runtime (bottom row) for the full experiment from \Cref{subsec:num_exp_mc} (left column) and a close up view of the first $\num[group-separator={,}]{2000}$ iterations (right column).
    Top: value of the QoI $f_\text{adapt}$ (blue dots) and evolution of the estimated expectation (orange) and variance (green) during the Monte Carlo run.
    Bottom: required time to evaluate the adaptive model $M_\text{adapt}$ for each $\mu$ chosen by the random sampling, and details on time spent for evaluation and building of the required submodels of $M_\text{adapt}$, compare \cref{alg:adaptive_output}.
}
	\label{fig:mc_results}
\end{figure}

The results of our numerical experiments\footnote{Performed on a computer equipped with a AMD Ryzen Threadripper 2950X CPU with 16 cores running at 3.50GHz each (and hyper-threading enabled) and 128GB of memory available.} %
are visualized in \cref{fig:mc_results}.
As we observe in \cref{fig:mc_results} (top right), even after $\num[group-separator={,}]{2000}$ samples, slight variations in the estimated variance can be observed, such that more samples are required to obtain accurate estimates.
Nevertheless, a relatively fast convergence of the estimators can be observed for the problem at hand.
However, the experiments also show that the adaptive algorithm is capable of producing fast and reliable results by frequently using the ML-ROM.
\par
In the beginning, the FOM was queried 12 times (blue bars in \cref{fig:mc_results}, bottom) until a sufficiently accurate RB-ROM could be built (orange).
Starting with the 13th evaluation, the RB-ROM was accurate enough, resulting in subsequent faster evaluations of the target quantity (yellow bars).
After $200$ ROM solves, the ML-ROM was trained for the first time, causing some optimization effort (dark green bars).
With this ML-ROM, subsequent evaluations are partly given by either the ML-ROM or the RB-ROM, depending on the certification of the respective model.
On average, a single solve using the FOM takes $216$s, while the RB-ROM requires $4.9$s on average and the ML-ROM $2.1$s.
After each additional $200$ ROM solves, i.e.\ $511$ and finally $1011$ solves in total, the ML-ROM is trained again, now using the larger available training data (dark green bars in \cref{fig:mc_results}, bottom).
The training times for the ML-ROM amount to $1231.1$s, $1256.6$s, and $3691.5$s.
Those additional training runs based on more data improve the ML-ROM and thus increase the ratio of ML-ROM vs.~RB-ROM evaluations:
The ratio of ML-ROM solves increases from $0\%$ before the first training to $35.9\%$ between the first and the second training.
After the second training, $59.4\%$ of the ML-ROM solutions are accepted, and finally, after the last training, for $97.6\%$ of the parameters the ML-ROM produces a solution that is accurate enough.
The total ML-ROM optimization time of $6179.3$s pays off after $2201$ uses of the ML-ROM-evaluation, which is quickly reached in view of the $97.6\%$ efficiency of the ML-ROM after the final, third optimization.
Note that the time to evaluate the FOM as well as the RB-ROM significantly increases for a larger number of time-steps, in contrast to the ML-ROM.
Thus, we expect even more pronounced speed-ups for experiments with larger $K\gg 1$.

\section{Conclusion}

We summarize by emphasizing two main contributions of the current presentation.
First, we have demonstrated a method for certification of ML-based models
for predicting input-output maps governed by parametric PDEs.
The main idea is to use as prediction target the coefficients of a reduced basis space.
This enables to apply the machinery of residual-based a posteriori certification
of ROMs to obtain a certified ML-model.
The guaranteed
error bound for each single prediction of ML-based models is a striking feature as
typically no strong accuracy statements can be given for learning-based models.
In particular, this allows to train with little data, as the prediction can be certified a posteriori.
Second, we have introduced and formalized a hierarchical and adaptive
RB-ML-ROM
and
demonstrated how it successfully can be implemented in different ways and be used
in many-query settings such as optimization as well as uncertainty quantification by Monte Carlo
simulation.
Especially the type of ML method to use is very flexible, we have demonstrated this by involving kernel methods as well
as neural networks.

Despite its generality, the abstract framework certainly has limitations.
The most limiting factor is that the success of the method still depends on the (usual) MOR assumption
that the parametric variation of the solutions is not too severe over the
parameters that are used as sample points.
In the extreme case, if there is no continuity of the parameter-to-solution
map $\mu \mapsto u(\mu)$, certainly no good speedup of the (hierarchical) RB-ML-ROM model can
be obtained, as the model will always refer to an evaluation of the FOM.
But remarkably even in this case, the model will still (trivially) be accurate. 
Also the RB-ML-ROM model has some limitations. Especially, in our particular case of
a reduced basis model
which is incrementally extended upon refinement, the main drawback clearly is 
that the reduced model is growing over time, hence the costs for RB-ROM evaluation as well as
the certification of both the ML and the RB sub-model by error estimation
will be increasing after every extension of the reduced basis.
Here, some ``unlearning'' of snapshots, or similar shrinking/subset
selection strategies could be devised.

As outlook, it might be interesting to investigate formal properties of the overall algorithm such as
convergence statements or even optimality of the adaptive model.
The model hierarchy might be extended by further layers
leading to even more sophisticated hierarchical models.

\end{document}